\documentclass[reqno]{amsart}

\usepackage{graphics}
\usepackage{amsthm}
\usepackage{amsmath}
\usepackage{amssymb}
\usepackage{stmaryrd}
\usepackage{mathrsfs}
\usepackage{enumerate}

\newcommand{\ud}[0]{\,\mathrm{d}}
\newcommand{\doo}[0]{\partial}

\newcommand{\intav}[0]{-\!\!\!\!\!\!\int}

\newcommand{\diverg}[0]{\operatorname{div}}
\newcommand{\dist}[0]{\operatorname{dist}}

\newcommand{\abs}[1]{\left|#1\right|}
\newcommand{\Babs}[1]{\Big|#1\Big|}
\newcommand{\norm}[2]{\left|#1\right|_{#2}}

\newcommand{\Norm}[2]{\left\|#1\right\|_{#2}}
\newcommand{\BNorm}[2]{\Big\|#1\Big\|_{#2}}
\newcommand{\bnorm}[2]{\Big|#1\Big|_{#2}}
\newcommand{\fNorm}[2]{\|#1\|_{#2}}
\newcommand{\pair}[2]{\left\langle #1,#2 \right\rangle}
\newcommand{\ave}[1]{\langle #1\rangle}

\newcommand{\bddlin}[0]{\mathscr{L}}
\newcommand{\kernel}[0]{\mathsf{N}}
\newcommand{\range}[0]{\mathsf{R}}
\newcommand{\domain}[0]{\mathsf{D}}

\newcommand{\local}[0]{\operatorname{loc}}
\newcommand{\BMO}[0]{\operatorname{BMO}}
\newcommand{\supp}[0]{\operatorname{supp}}
\newcommand{\Car}[0]{\operatorname{Car}}

\renewcommand{\Re}[0]{\operatorname{Re}}

\newcommand\R{\mathbf{R}}
\newcommand\C{\mathbf{C}}
\newcommand\N{\mathbf{N}}
\newcommand\Z{\mathbf{Z}}

\newcommand{\Exp}[0]{\mathbb{E}}

\newcommand{\radem}[0]{\varepsilon}
\newcommand{\rbound}[0]{\mathscr{R}}
\newcommand{\Rad}[0]{\operatorname{Rad}}

\theoremstyle{plain}
\newtheorem{theorem}{Theorem}[section]
\newtheorem{proposition}[theorem]{Proposition}
\newtheorem{corollary}[theorem]{Corollary}
\newtheorem{lemma}[theorem]{Lemma}

\theoremstyle{definition}
\newtheorem{definition}[theorem]{Definition}

\theoremstyle{remark}
\newtheorem{remark}[theorem]{Remark}

\begin{document}

\title{Kato's square root problem in Banach spaces}

\author[T.~Hyt\"onen]{Tuomas Hyt\"onen}
\address{Department of Mathematics and Statistics, University of Helsinki, Gustaf H\"allstr\"omin katu 2b, FI-00014 Helsinki, Finland}
\email{tuomas.hytonen@helsinki.fi}

\author[A.~McIntosh]{Alan McIntosh}
\address{CMA, Australian National University, Canberra ACT 0200, Australia}
\email{Alan.McIntosh@maths.anu.edu.au}

\author[P.~Portal]{Pierre Portal}
\address{CMA, Australian National University, Canberra ACT 0200, Australia}
\email{pierre.portal@maths.anu.edu.au}

\date{26 February 2007.}


\begin{abstract}
Let \(L\) be an  elliptic differential operator with bounded measurable coefficients, acting in Bochner spaces \(L^{p}(\R^{n};X)\) of \(X\)-valued functions on \(\R^n\). We characterize Kato's square root estimates \(\|\sqrt{L}u\|_{p} \eqsim \|\nabla u\|_{p}\) and the \(H^{\infty}\)-functional calculus of~\(L\) in terms of \(R\)-boundedness properties of the resolvent of~\(L\), when \(X\) is a Banach function lattice with the UMD property, or a noncommutative \(L^{p}\) space. To do so, we develop various vector-valued analogues of classical objects in Harmonic Analysis, including a maximal function for Bochner spaces. In the special case \(X=\C\), we get a new approach to the \(L^p\) theory of square roots of elliptic operators, as well as an \(L^{p}\) version of Carleson's inequality.
\end{abstract}

\maketitle

\tableofcontents

\section{Introduction}\label{intro}

The development of a theory of singular integrals for vector functions, which take their values in an infinite-dimensional Banach space, may be viewed as an accelerated replay --- with new actors, insight, and considerable improvisation --- of the original development in the scalar-valued setting. During the 1980's, this theory advanced from D.~L.\ Burkholder's \cite{Burkholder:83} extension of M.~Riesz' classical theorem on the Hilbert transform boundedness, via J.~Bourgain's \cite{Bourgain:86}, T.~R.\ McConnell's \cite{McConnell:84} and F.~Zimmermann's \cite{Zimmermann} results on Calder\'on--Zygmund principal value convolutions and Marcinkiewicz--Mihlin multipliers, to T.~Figiel's \cite{Figiel:90} vector-valued generalization of the  $T(1)$ theorem of G.~David and J.-L.\ Journ\'e. More recently, there has been a new boom of activity in developing the vector-valued estimates to match the needs of a wide variety of applications especially in the field of Partial Differential Equations. An important opening move into this direction was made by L.~Weis \cite{Weis}; further developments and references are recorded in~\cite{DHP,kunstweis}.

The aim of the present paper is to continue the vector-valued program so as to catch up with some of the latest achievements in scalar-valued Harmonic Analysis. More precisely, we are going to develop a Banach space theory for the square roots of elliptic operators appearing in the famous problem of T.~Kato, which was recently solved by P.~Auscher, S.~Hofmann, M.~Lacey, A.~McIntosh and Ph.~Tchamitchian \cite{AHLMT}, and more generally for the perturbed Dirac operators treated in a subsequent work by A.~Axelsson, S.~Keith and A.~McIntosh \cite{AKM}. These objects are no longer Calder\'on--Zygmund operators, and may even fail to have a pointwise defined kernel.

For this reason, their study is considered a move beyond Calder\'on--Zygmund theory. In the scalar valued \(L^{p}\) case, this has recently attracted much attention. An extrapolation technique developed by S.~Blunck and P.~Kunstmann \cite{blunckkunst} allows to extend \(L^{2}\) results to the \(L^{p}\) setting for \(p\) in an open interval \((p_{-},p_{+})\), which may be strictly smaller than the whole reflexive range $(1,\infty)$ admissible for classical operators. P.~Auscher's memoir \cite{auscher} presents the large range of applications of this method and demonstrates that the \(L^{p}\) behavior of objects associated with an elliptic operator \(L\) (its functional calculus, Riesz transforms, square functions, etc.) is ruled by four critical numbers: \(p_{-}(L),p_{+}(L)\) (the limits of the range of $p$'s for which the semigroup \((e^{-tL})_{t>0}\) is $L^p$-bounded), and
 \(q_{-}(L),q_{+}(L)\) (the limits of the range of $p$'s for which the family \((\sqrt{t}\nabla e^{-tL})_{t>0}\) is $L^p$-bounded).
In a recent series of papers by P.~Auscher and J.~M.\ Martell \cite{auschermartell}, these results are extended to a more general setting, allowing weighted estimates on spaces of homogeneous type. We also refer to their papers for the history of these developments.

Our work takes a different approach. Since we are aiming at a Banach space-valued theory, where no easier $L^2$ case is available as a starting point, we cannot rely on an extrapolation method, but need to work directly in the spaces $L^p(\R^n;X)$.  It is interesting, even in the scalar case $X=\C$, to see that the methods from \cite{AHLMT} and \cite{AKM} can in fact be extended to an \(L^{p}\) situation.
This requires a set of new techniques. We develop, in particular, a Banach space valued analogue of the ``reduction to the principal part'' method used to solve Kato's problem (Theorem \ref{reduction}). This is based on adequate off-diagonal estimates (Proposition \ref{offdiag}), and on the fact that resolvents of an unperturbed Hodge-Dirac operator are, in some sense, equivalent to conditional expectations with respect to the dyadic filtration of \(\R^n\) (Corollary \ref{PisA}).
This result, which is handled in the classical case by a $T(1)$ Theorem for Carleson measures (see \cite{at}), is obtained in our context by extending ideas from \cite{AKM}. To do so, we develop Banach space valued analogues of classical estimates such as Poincar\'e's inequality, and Schur's Lemma.

Finally we establish an analogue of Carleson's inequality (Theorem~\ref{carlesonestimate}) to handle the principal part.
This is a crucial step and requires the \(L^{p}\) boundedness of an appropriate (Rademacher) maximal function which we introduce and study in Section \ref{maximal}. We prove its boundedness in \(L^{p}(\R^{n};X)\) when $1<p<\infty$ provided that \(X\) is either a UMD function lattice, or a non commutative \(L^{q}\) space for some \(1<q<\infty\), or a space with Rademacher type~\(2\). We thus obtain a satisfying result in most of the concrete spaces of interest, but the boundedness of the Rademacher maximal function (and hence the Kato estimates) in general UMD valued Bochner spaces remains open.

The paper is organized as follows. 
In Section~\ref{prelim} we provide the reader with a concise introduction to the concepts and results from the theory of Banach spaces and Banach space valued Harmonic Analysis used in this paper. Section~\ref{result} contains the statements of the main results, and their reduction to the main estimate which is then dealt with in the rest of the paper. We develop vector-valued analogues of various classical results, which came to use in the proof of the scalar Kato problem, in Section~\ref{sec:props}. Section~\ref{unpert} deals with the Banach space valued analogues of classical inequalities associated with an unperturbed Hodge--Dirac operator, and in particular the relationship with the dyadic conditional expectations. In Section~\ref{sect:reduction} we reduce the main estimate to its principal part. Our Rademacher maximal function is studied in Section~\ref{maximal} and applied in Section~\ref{carleson} to prove an analogue of Carleson's inequality. This is used to reduce the principal part estimate to an analogue of a Carleson measure condition, which is finally verified in Section~\ref{sec:conclusion} by essentially the same stopping time argument as in~\cite{AHLMT} and~\cite{AKM}.

Additional results are presented in three appendices. In Appendix~\ref{app:IplusK} we show how the assumptions of the main theorem can in some cases be checked under appropriate ellipticity conditions. In Appendix~\ref{app:paraproducts} we relate our Carleson inequality to the boundedness of vector-valued paraproducts, and finally Appendix~\ref{app:l1} contains a counterexample related to the Rademacher maximal function.


\section{Preliminaries}\label{prelim}

This work is concerned with resolvent bounds, $H^{\infty}$ functional calculus, and quadratic estimates for certain partial differential operators acting in $L^p$ spaces of Banach space valued functions. In order to streamline the actual discussion, we start by recalling the relevant notions and a number of results which will be repeatedly used in the sequel.

To express the typical inequalities ``up to a constant'' we use the notation \(a \lesssim b\) to mean that there exists \(C<\infty\) such that \(a \leq Cb\), and the notation \(a \eqsim b\) to mean that \(a \lesssim b \lesssim a\). The implicit constants are meant to be independent of other relevant quantities. If we want to mention that the constant \(C\) depends on a parameter \(p\), we write \(a \lesssim _{p} b\).

\begin{definition}
Let \(A\) be a closed operator acting in a Banach space \(Y\). It is called \emph{bisectorial} with angle \(\theta\) if its spectrum \(\sigma(A)\) is included in a bisector:
\begin{equation*}\begin{split}
  \sigma(A)\subseteq S_{\theta} &:= \Sigma_{\theta} \cup\{0\}\cup (-\Sigma_{\theta}),\quad\text{where}\\
 \Sigma_{\theta} &:= \{z \in \C\setminus\{0\} \;;\; |\arg(z)| < \theta\},
\end{split}\end{equation*}
and outside the bisector it verifies the following resolvent bounds:
\begin{equation}\label{eq:biSect}
  \forall \theta' \in (\theta,\frac{\pi}{2}) \quad \exists C>0 \quad \forall \lambda \in \C\setminus S_{\theta'} \quad
  \Norm{\lambda(\lambda I-A)^{-1}}{\bddlin(Y)} \leq C.
\end{equation}
\end{definition}
We often omit the angle, and say that \(A\) is bisectorial if it is bisectorial with \emph{some} angle \(\theta \in [0,\frac{\pi}{2})\). 
One sees that \(A\) is bisectorial if and only if it satisfies the resolvent bound in \eqref{eq:biSect} on the imaginary axis, i.e.,
\begin{equation*}
  \|(I+itA)^{-1}\|\leq C,\qquad t\in\R.
\end{equation*}

For \(0<\nu<\pi/2\), let \(H^{\infty}(S_{\nu})\) be the space of bounded functions on \(S_{\nu}\), which are holomorphic in \(S_{\nu}\setminus\{0\}\), and consider the following subspace of functions with decay at zero and infinity:
\begin{equation*}\begin{split}
  H_{0} ^{\infty}(S_{\nu}) := \Big\{&\phi \in H^{\infty}(S_{\nu}):\\
     &\exists \alpha,C\in(0,\infty)\quad \forall z \in S_{\nu} \quad \abs{\phi(z)} \leq C\abs{\frac{z}{1+z^{2}}}^{\alpha} \Big\}.
\end{split}\end{equation*}
For a bisectorial operator \(A\) with angle \(\theta<\omega<\nu<\pi/2\), and \(\psi\in H^{\infty}_0(S_{\nu})\), we define
\begin{equation*}
  \psi(A)u := \frac{1}{2i\pi} \int_{\partial S_{\omega}} \psi(\lambda)(\lambda-A)^{-1}u \ud\lambda,
\end{equation*}
where \(\partial S_{\omega}\) is parameterized by arclength and directed anti-clockwise around \(S_{\omega}\).

\begin{definition}
Let \(A\) be a bisectorial operator with angle \(\theta\), and \(\nu \in (\theta,\frac{\pi}{2})\). \(A\) is said to admit a {\em bounded \(H^{\infty}\) functional calculus with angle~\(\nu\)} if
\(
\quad \exists C<\infty \quad \forall \psi \in H_{0} ^{\infty}(S_{\nu}) \quad \|\psi(A)y\|_{Y} \leq C\|\psi\|_{\infty} \|y\|_{Y}.
\)
\end{definition}
On the closure \(\overline{\range(A)}\) of the range space \(\range(A)\), we then define a bounded operator \(f(A)\), for every \(f \in H^{\infty}(S_{\nu})\), by
\(f(A)u=\underset{n \to \infty}{\lim}\psi_{n}(A)u\), where \(\psi_{n} \in H_{0} ^{\infty}(S_{\omega})\) are uniformly bounded and tend to \(f\) locally uniformly on \(S_{\omega}\setminus\{0\}\).
In a reflexive Banach space, there holds \(X=\kernel(A)\oplus\overline{\range(A)}\) (cf.~\cite{Haase}, Proposition~2.1.1, for the sectorial case which is readily adapted to the present context), so that denoting by \(\mathbb{P}^{0}\) the associated projection onto the null space \(\kernel(A)\), we can finally define the bounded operator \(f(A)\) by
\begin{equation*}
  f(A)u = f(0)\mathbb{P}^{0}u+\underset{n \to \infty}{\lim}\psi_{n}(A)u.
\end{equation*}

We also often omit the angle and just say that \(A\) has an  \(H^{\infty}\) functional calculus.  The detailed construction of this calculus, and much more information, can be found in \cite{CDMY,Haase,kunstweis}.

A crucial aspect of the functional calculus is its harmonic analytic characterization.
If \(Y\) is a Hilbert space, it is shown in \cite{mc86} that \(A\) has an \(H^{\infty}\) functional calculus with angle \(\nu\) if and only if the following  {\em quadratic estimate} holds 
\begin{equation*}
  \Big(\int_{0} ^{\infty} \Norm{\psi(tA)y}{Y}^{2} \frac{\ud t}{t}\Big)^{1/2} \eqsim \Norm{y}{Y}
\end{equation*}
{\em for some} non-zero function \(\psi \in H_{0} ^{\infty} (S_{\nu})\).
In the space \(L^{p}(\R^{n};\C)\) (\(1<p<\infty\)), it has been shown in \cite{CDMY} that the above norms need to be replaced by 
\begin{equation*}
  \BNorm{\Big(\int_{0} ^{\infty} \abs{\psi(tA)y}^{2} \frac{\ud t}{t}\Big)^{1/2}}{p}
\end{equation*}
as in the Littlewood--Paley theory. In a general Banach space, the correct characterization involves randomized sums of the form
\begin{equation*}
  \Exp\BNorm{\sum_{k\in\Z} \varepsilon_{k} \psi(2^{k}A)y}{Y},
\end{equation*}
where \((\varepsilon_{k})_{k \in \Z}\) are independent Rademacher variables on some probability space $\Omega$ (i.e., they take each of the two values $+1$ and $-1$ with probability $1/2$), and $\Exp$ is the mathematical expectation. 
These randomized norms provide the right analogue of the quadratic norms used in \(L^{p}\) and for this reason, somewhat loosely speaking, we will occasionally also refer to inequalities for the randomized norms as ``quadratic estimates''.

\begin{proposition}[Khintchine--Kahane inequalities]
\label{kk-inequalities}
Let \(Y\) be a Banach space, and \((y_{k})_{k \in \Z} \subset Y\). Then for each \(1<p<\infty\), there exists \(C_{p}>0\) such that
\begin{equation*}
  \mathbb{E}\BNorm{\sum_k \varepsilon_{k} y_{k}}{Y}
  \leq \Big(\mathbb{E}\BNorm{\sum_k \varepsilon_{k} y_{k}}{Y}^{p}\Big)^{1/p}
  \leq C_p\mathbb{E}\BNorm{\sum_k \varepsilon_{k} y_{k}}{Y}.
\end{equation*}
Moreover, if \(Y =  L^{q}\) for some \(1<q<\infty\) (or more generally a Banach lattice with finite cotype), then
\begin{equation*}
  \mathbb{E}\BNorm{\sum_k \varepsilon_{k} y_{k}}{Y} \eqsim
   \BNorm{\Big(\sum_k |y_{k}|^{2}\Big)^{1/2}}{Y}.
\end{equation*}
\end{proposition}

When using such randomized sums, it is often convenient to introduce the space \(\Rad(Y)\) of sequences \((y_{k})_{k \in \Z} \subset Y\) such that
\(\sum_{|k|<n} \varepsilon_{k} y_{k}\) converges in \(L^1(\Omega;Y)\), with the norm defined by 
\begin{equation*}
  \Norm{(y_{k})_{k \in \Z}}{\Rad(Y)} = \mathbb{E}\BNorm{\sum_k \varepsilon_{k} y_{k}}{Y}.
\end{equation*}

These norms involve discrete rather than continuous sums, but this techincal difference is unimportant. In fact, we could avoid discretization by using Banach space valued stochastic integrals as in \cite{tuomasLPS}, but this would only add an unnecessary level of complexity. 
An important problem, however, is the fact that the quadratic norms are not, outside the Hilbertian setting, independent of the choice of \(\phi \in H_{0} ^{\infty} (S_{\theta})\). To ensure such an independence, one has to assume (see \cite{kunstweis}) that the family
\(\{\lambda(\lambda I-A)^{-1} \;;\; \lambda \not \in S_{\theta}\}\) is not only bounded (bisectoriality) but \(R\)-bounded (\(R\)-bisectoriality) in the following sense.

\begin{definition}
Let \(X\) be a Banach space. A family of bounded linear operators \(\Psi \subset \mathcal{L}(X)\) is called {\em
R-bounded} if there exists a constant \(C\) such that for all\(N \in \N\), \(T_{1},...,T_{N} \in \Psi\), and \(x_{1},...,x_{N}
\in X\), there holds
\begin{equation*}
  \Exp \BNorm{\sum_{j=1} ^{N}
    \radem_{j}T_{j}x_{j}}{}
  \leq C \Exp \BNorm{\sum_{j=1} ^{N} \radem_{j}x_{j}}{}.
\end{equation*}
\end{definition}

A uniformly bounded family of operators is not necessarily R-bounded, as can be seen by considering translations on \(L^{p}\), \(p\neq 2\). In fact, the property that every uniformly bounded family is R-bounded characterizes Hilbert spaces  up to isomorphism. 
This is in contrast to the scalar multiplication where Kahane's principle holds:

\begin{proposition}[Contraction principle]
\label{contraction-principle}
Let \(X\) be a Banach space, and \(\lambda=(\lambda_{k})_{k \in \Z} \in \ell^{\infty}\). 
Then \(\forall N \in \N, \quad \forall x_{1},...,x_{N} \in X\)
\begin{equation*}
  \mathbb{E} \Big\| \sum \limits _{j=1} ^{N}
  \varepsilon_{j}\lambda_{j}x_{j}\Big\| \leq 2 \|\lambda\|_{\infty}
  \mathbb{E} \Big\| \sum
  \limits _{j=1} ^{N} \varepsilon_{j}x_{j}\Big\|.
\end{equation*}
\end{proposition}

An immediate but useful consequence of Propositions \ref{contraction-principle} and \ref{kk-inequalities} is the following (see e.g.~\cite{kunstweis}).

\begin{proposition}
Let \(X\) be a Banach space, and \((f_{k})_{k \in \Z}\subset L^{\infty}(\R^{n})\) be a bounded sequence of functions.
Then the family of multiplication operators defined by \(T_{k}u = f_{k}u\) is R-bounded on \(L^{p}(\R^{n};X)\) for all \(1<p<\infty\).
\end{proposition}

The concept of \(R\)-boundedness is crucial in Banach space valued Harmonic Analysis. It is described in detail in \cite{kunstweis}, where the following characterization can also be found (see Section 12 of \cite{kunstweis}):

\begin{theorem}[Kalton, Kunstmann, Weis]
Let \(Y\) be a UMD Banach space, and \(A\) be an R-bisectorial operator acting on \(Y\). 
Then \(A\) has an \(H^{\infty}\) functional calculus if and only if
\begin{equation*}
\begin{cases}
\underset{1 \leq |t| \leq 2}{\sup}
  \Exp\BNorm{\sum \limits _{k \in \Z} \varepsilon_{k}2^{k}tA(I+(2^{k}tA)^{2})^{-1}y}{Y} \lesssim \|y\|_{Y} \quad
\forall y \in Y, \\
\underset{1 \leq |t| \leq 2}{\sup} \Exp\BNorm{\sum \limits _{k \in \Z}
   \varepsilon_{k}2^{k}tA^{*}(I+(2^{k}tA^{*})^{2})^{-1}y^{*}}{Y^{*}} \lesssim \|y^{*}\|_{Y^{*}} \quad
\forall y^{*} \in Y^{*}. \\
\end{cases}
\end{equation*}
\end{theorem}

The main body of this paper is concerned with proving this kind of estimates when \(Y=L^p(\R^n;X^N)\) is the Bochner space of functions with values in the Cartesian product \(X^N\) of \(N\) copies of a Banach space \(X\), and \(A\) is a \emph{perturbed Hodge--Dirac} operator, as defined in the next section. Let us only mention at this point that our operators will be the ``simplest'' extensions of the classical Hodge--Dirac operators to the Banach space valued setting, namely tensor products \(T \otimes I_{X}\) of an operator \(T\) acting in \(L^p(\R^n;\C^N)\) with the identity \(I_X\). The study of such operators is by no means trivial. Already in the case when \(T\) is the possibly simplest singular integral operator, the Hilbert transform, the boundedness of \(T\otimes I_{X}\) in \(L^{p}(\R;X)\) is equivalent to \(X\) being a so-called \emph{UMD space}, which means the \emph{u}nconditional convergence of \emph{m}artingale \emph{d}ifference sequences in \(L^p(\Omega;X)\) for \(1<p<\infty\) and \(\Omega\) any probability space.

This class of spaces  is the most important one for vector-valued Harmonic Analysis. All UMD spaces are reflexive (and even super-reflexive; cf.~\cite{Bourgain:83}). The principal examples include the reflexive Lebesgue, Lorentz, Sobolev, and Orlicz spaces, as well as the reflexive \emph{noncommutative \(L^p\) spaces}. A recent survey paper on UMD spaces is \cite{Burkholder:UMD}. The abovementioned equivalence with the Hilbert transform boundedness, due in one direction to Burkholder \cite{Burkholder:83} and in the other to Bourgain \cite{Bourgain:83}, lies at the heart of the theory, and is characteristic of the interaction between probabilistic and analytic methods. It is, for instance, needed in the proof of the following multiplier theorem, which we often resort to in the sequel. The original statement of this kind was obtained by  Bourgain~\cite{Bourgain:86} and McConnell~\cite{McConnell:84}, but the somewhat more general formulation given here is due to Zimmermann~\cite{Zimmermann}. 

\begin{theorem}[Bourgain, McConnell, Zimmermann]
\label{basicmult}
Let $n\geq 1$. If \textup{(}and only if\textup{)} $X$ is a UMD
space and $1<p<\infty$, then every symbol
$m:\R^n\setminus\{0\}\to \C$ such that
\begin{equation*}
  \sup\{\abs{\xi}^{\abs{\alpha}}
    D^{\alpha}m(\xi):\alpha\in\{0,1\}^n,
    \xi\in\R^n\setminus\{0\}\}<\infty
\end{equation*}
gives rise to a bounded Fourier multiplier $T_m\in\bddlin(L^p(\R^n,X))$ defined by 
$\mathcal{F}(T_m u)(\xi) = m(\xi) \mathcal{F}(u)(\xi)$, where \(\mathcal{F}\) denotes the Fourier transform.
\end{theorem}

With somewhat stronger conditions on the symbol, we also have stronger conclusions. Let us say that a symbol $m:\R^n\to \C$ has \emph{bounded variation} if for some $C<\infty$ and all $\alpha\in\{0,1\}^n$, there holds
\begin{equation*}
  \int_{\R}\cdots\int_{\R} \abs{D^{\alpha}m(\xi)}\ud\xi^{\alpha}\leq C<\infty,
\end{equation*}
where the integration is with respect to all the variables $\xi_i$ such that $\alpha_i=1$, and the estimate is required uniformly in the remaining variables $\xi_j$. (The case $\alpha=0$ is understood as the boundedness of $m(\xi)$ by $C$.) We say that a collection of symbols $\mathscr{M}$ has uniformly bounded variation if the symbols $m\in\mathscr{M}$ satisfy this condition with the same $C$. See \cite{kunstweis} for the proof of the following useful result:

\begin{proposition}\label{RbddMult}
Let $n\geq 1$, $X$ be a UMD, and $1<p<\infty$. Let $\mathscr{M}$ be a collection of symbols of uniformly bounded variation.
Then the collection of Fourier multipliers $T_m$, $m\in\mathscr{M}$, is an $R$-bounded subset of $\bddlin(L^p(\R^n;X))$.
\end{proposition}

Another important estimate in UMD spaces, analogous to the previous one, is the following $R$-boundedness of conditional expectations.  It is an extension of a classical quadratic estimate due to Stein~\cite{Stein}, which was found in the vector-valued situation by Bourgain \cite{Bourgain:86}. See also \cite{FigWoj} for a proof.

\begin{proposition}[Stein's inequality]
\label{rbdd-At}
Let \(X\) be a UMD Banach space, \((\Omega, \Sigma, \mu)\) a measure space, and \(1<p<\infty\). 
Then any  increasing  sequence of conditional expectations on \(L^{p}(\Omega;X)\) is R-bounded.
\end{proposition}

We will mostly be concerned with the conditional expectations related to the dyadic filtration of $\R^n$. This is defined by the system of \emph{dyadic cubes}
\begin{equation*}
  \triangle = \bigcup_{k \in \Z}\triangle_{2^k},\qquad
  \triangle_{2^k}:=\big\{2^k([0,1)^n+m):m\in\Z^n\big\}.
\end{equation*}
The corresponding conditional expectation projections are denoted by
\begin{equation*}
  A_{2^k} u(x):=\ave{u}_Q:=\intav_Q u(y)\ud y:=\frac{1}{\abs{Q}}\int_Q u(y)\ud y,\qquad
  x\in Q\in\triangle_{2^k}.
\end{equation*}
The integral average notation above will also be used with other measurable sets from time to time.

Other important Banach space properties are the following:

\begin{definition}
Let \(X\) be a Banach space, and \(1\leq t\leq 2\leq s\leq\infty\). Then \(X\) is said to have \emph{(Rademacher) type} \(t\) if
\begin{equation*}
  \Exp\BNorm{\sum_{k\in\Z} \radem_k x_k}{X}\lesssim
  \Big(\sum_{k\in\Z}\Norm{x_k}{X}^t\Big)^{1/t}
\end{equation*}
for all \(x_k\in X\), and \emph{(Rademacher) cotype} \(s\) if
\begin{equation*}
  \Big(\sum_{k\in\Z}\Norm{x_k}{X}^s\Big)^{1/s}\lesssim
  \Exp\BNorm{\sum_{k\in\Z} \radem_k x_k}{X}
\end{equation*}
for all \(x_k\in X\), where the usual modification is understood if \(s=\infty\). The space is said to have \emph{nontrivial type} if it has some type $t>1$, and \emph{nontrivial}, or \emph{finite}, \emph{cotype} if it has some cotype \(s<\infty\).
\end{definition}

These conditions become stronger with increasing $t$ and decreasing $s$, and only Hilbert spaces (up to isomorphism) enjoy both the optimal type and cotype $t=s=2$.
For the present purposes, the most important thing is to know that every UMD space has both nontrivial type and cotype. The property of finite cotype is also characterized (see \cite{DJT}, 12.27) by the comparability of Rademacher and Gaussian random sums,
\begin{equation}\label{eq:RadVsGauss}
  \Exp\BNorm{\sum_{k\in\Z}\radem_k x_k}{X}
  \eqsim\Exp\BNorm{\sum_{k\in\Z}\gamma_k x_k}{X}\quad\Leftrightarrow\quad
  X\text{ has finite cotype},
\end{equation}
where the \(\gamma_k\) are independent random variables with the standard normal distribution.

These notions, as well as the Khintchine--Kahane inequalities \ref{kk-inequalities}, are central in a circle of ideas which can be roughly referred to as ``averaging in Banach spaces'', and which forms the core of vector-valued harmonic analysis.
A gentle introduction to this topic can be found in \cite{zebook}.

In addition to the above conditions, which are well known in the theory of Banach spaces, we need to introduce a new class of spaces, the defining property of which is the boundedness of the following \emph{Rademacher maximal function}:
\begin{equation*}\begin{split}
  M_Ru(x):=\sup\Big\{\Exp &\BNorm{\sum_{k\in\Z}\radem_k\lambda_k A_{2^k}u(x)}{X}: \\
    &\lambda=(\lambda_k)_{k\in\Z}\text{
    finitely non-zero with }\Norm{\lambda}{\ell^2(\Z)}\leq 1\Big\}.
\end{split}\end{equation*}
Note that, under the identification $X\eqsim\bddlin(\C,X)$, this is the
$R$-bound of the set
\begin{equation*}
  \big\{A_{2^k}u(x):k\in\Z\big\}
  =\big\{\ave{u}_Q:Q\owns x\big\}.
\end{equation*}
In particular, if $X$ is a Hilbert space, we recover the usual dyadic maximal function.

\begin{definition}
We say that the Banach space $X$ has the RMF property, if \(M_{R}\) is bounded from \(L^{2}(\R^{n};X)\) to \(L^{2}(\R^{n})\).
\end{definition}

We do not yet completely understand  how this new class of spaces relates to the other Banach space notions discussed above, which forces us to adopt this property as an additional assumption. It would be particularly useful to know if every UMD space has RMF, since this would allow us to state our main theorem in the generality of all UMD spaces, but the question remains open. However, in Section~\ref{maximal} we show that the RMF property does hold in most of the concrete situations of interest. The classes of Banach spaces appearing in the statement are also defined in Section~\ref{maximal}. 

\begin{proposition}\label{radmax}
A Banach space which is a UMD function lattice, or a noncommutative \(L^{p}\) space for \(1<p<\infty\), or which has Rademacher type \(2\), has RMF.
\end{proposition}


\section{Statement of the results}
\label{result}

The square root problem originally posed by T.~Kato was an
operator-theoretic question in an abstract Hilbert space, but it was
observed  in \cite{lions} and \cite{mc}  that the desired estimate
was invalid in this generality (see \cite{AHLMT} for references and more historical information). This shifted
the attention towards more concrete differentiation and
multiplication operators in $L^2(\R^n;\C^N)$, ones of interest in
the actual applications that Kato had in mind when formulating his
problem. Our Banach space framework is obtained by modifying the concrete Kato problem, so as the replace $\C^N$ by $X^N$, and $L^2$ by $L^p$. The various differentiation and multiplication operators are simply replaced by their natural tensor extensions acting on $X$-valued functions. The set-up, which we now present in detail, is closely related to that of \cite{AKM}, Section~3.

Let $X$ be a Banach space, $1<p<\infty$, and $n,\ n_1,\ n_2,\ N\in\Z_+$ with $N=n_1+n_2$. Let $D$ be a homogeneous first order partial differential operator with constant $\bddlin(\C^{n_1},\C^{n_2})$-coefficients, and $D^*$ be its adjoint. We assume that
\begin{equation}\label{eq:DDstar}
  DD^*D=-\Delta D.
\end{equation}
The principal case of interest is
\[
  \{n_1,n_2,D,D^*\}
  =\{1,n,\nabla,-\diverg\},
\]
but it is convenient to consider the abstract formulation, because it makes the assumptions symmetric in $D$ and $D^*$. (Note that \eqref{eq:DDstar} is equivalent to the similar equation with $D$ and $D^*$ reversed by taking adjoints of both sides.) For $i=1,2$, let $A_i\in L^{\infty}(\R^n;\bddlin(\C^{n_i}))$ be bounded matrix-valued functions, which we identify with multiplication operators on $L^p(\R^n;X^{n_i})$ in the natural way.  We assume the estimate
\begin{equation*}
  \Norm{A_i}{L^{\infty}(\R^n;\bddlin(\C^{n_i}))}
  +\Norm{A_i^{-1}}{L^{\infty}(\R^n;\bddlin(\C^{n_i}))}
 \leq C,\qquad i=1,2.
\end{equation*}

In the space
\begin{equation*}
  L^p(\R^n;X^N)\equiv
  L^p(\R^n;X^{n_1})\oplus L^p(\R^n;X^{n_2})
\]
 we consider the operators
\begin{equation*}
  \Gamma=\begin{pmatrix}
    0 & 0 \\ D & 0 \end{pmatrix},\quad
  \Gamma^*=\begin{pmatrix}
    0 & D^* \\ 0 & 0 \end{pmatrix},\quad
  B_1=\begin{pmatrix}
    A_1 & 0 \\ 0 & 0 \end{pmatrix},\quad
  B_2=\begin{pmatrix}
    0 & 0 \\ 0 & A_2 \end{pmatrix}.
\end{equation*}
The first two are closed and nilpotent (i.e., the range \(\range(\Gamma)\subseteq\kernel(\Gamma)\), the null space; and the same with \(\Gamma^*\)) operators with their natural dense domains \(\domain(\Gamma)\) and \(\domain(\Gamma^*)\), while the latter two are everywhere defined and bounded.

The sum
\begin{equation*}
  \Pi = \Gamma + \Gamma^{*}
  =\begin{pmatrix} 0 & D^* \\ D & 0 \end{pmatrix}
\end{equation*}
is called the {\em Hodge-Dirac} operator.
Modified sums  of the form
\begin{equation*}\begin{split}
  \Pi_{B} &= \Gamma+\Gamma_B^*
   = \Gamma + B_{1}\Gamma^{*}B_{2}, \\
  \Pi_{B^{*}} &= \Gamma^* + \Gamma_{B^*}
   = \Gamma^{*} + B_{2}\Gamma B_{1}
\end{split}\end{equation*}
are then called {\em perturbed Hodge-Dirac} operators. It follows from general Operator Theory, using only the closedness or boundedness of the appropriate operators and the form of the matrices, that $\Pi_B$ and $\Pi_{B^*}$ are also closed and densely defined.

In the Hilbert space setting of \cite{AKM}, appropriate ellipticity conditions on \(B_{1}\) and \(B_{2}\) further imply, still by abstract operator theoretic methods, the defining resolvent estimates for the ($R$-)bisectoriality of \(\Pi_B\) and \(\Pi_{B^*}\). In the present situation, this is no longer the case; in fact, already when $X=\C$ but $p\neq 2$,  there exist elliptic second order differential operators which are not sectorial in \(L^{p}(\R^{n};\C)\) for some values of \(p\) (see \cite{ACT}). Thus we need to redefine the problem slightly, so as to adopt the analogues of some of the operator-theoretic conclusions in \cite{AKM} as the assumptions for our Harmonic Analysis. In particular, we \emph{assume} the existence of the following resolvents of $\Pi_B$ for all $t\in\R$:
\begin{equation}\label{eq:resolvents}\begin{split}
  R_t^B &:=(I+it\Pi_B)^{-1},\\
  P_t^B &:=(I+t^2\Pi_B^2)^{-1}=\frac{1}{2}(R_t^B+R_{-t}^B)
   =R_t^B R_{-t}^B,\\
  Q_t^B &:=t\Pi_B P_t^B=t\Pi_B(I+t^2\Pi_B^2)^{-1}
   =\frac{i}{2}(R_t^B-R_{-t}^B).
\end{split}\end{equation}

We can now state our main result.

\begin{theorem}\label{mainthm}
Let \(X\) be a UMD Banach space such that both $X$ and $X^*$ have RMF.
Let $\Pi_B$ and $\Pi_{B^*}$ be perturbed Hodge-Dirac operators defined in $L^p(\R^n;X^N)$ for all $p\in(p_-,p_+)\subseteq(1,\infty)$. Then the following are equivalent:
\begin{equation}\label{eq:PiBRsect}
  \Pi_B,\ \Pi_{B^*}\text{ are $R$-bisectorial in }
  L^p(\R^n;X^N) \text{ for all } p \in (p_{-},p_{+}).
\end{equation}
\begin{equation}\label{eq:PiBfc}
  \Pi_B,\ \Pi_{B^*}\text{ have $H^{\infty}$-calculus in }
  L^p(\R^n;X^N)\text{ for all } p \in (p_{-},p_{+}).
\end{equation}
\end{theorem}

The reason why we are forced to formulate this theorem for \(L^p\) estimates valid on open intervals of exponents, instead of  an individual \(p\), comes from the limitations in one particular step of the proof (our \(L^p\) version of Carleson's inequality); this will be discussed in somewhat more detail in Section~\ref{sec:conclusion}. Note that we do not require that \(2\in(p_-,p_+)\) here, whereas this is often the case in the scalar-valued results which are based on extrapolation of the \(L^2\) estimates.

The next corollary makes the relation to the square roots of second-order differential operators more explicit.

\begin{corollary}\label{maincor}
Let \(X\) be a UMD Banach space such that both $X$ and $X^*$ have RMF.
Let \(A\) and \(A^{-1}\) be multiplications by \(L^{\infty}(\R^{n};\bddlin(\C^n))\) functions, and \(L=-\diverg A\nabla\) be a sectorial operator in \(L^p(\R^n;X)\) for all \(p\in(p_{-},p_{+})\subseteq (1,\infty)\).
Then the following are equivalent:
\begin{equation}\label{eq:fourRbdd}\left\{\begin{split}
   &\text{The sets } \{(I+t^2L)^{-1}\}_{t>0},\ \{t\sqrt{-\Delta}(I+t^2L)^{-1}\}_{t>0},\\
   &\{(I+t^2L)^{-1}t\sqrt{-\Delta}\}_{t>0}
    \text{ and }\{t\sqrt{-\Delta}(I+t^2L)^{-1}t\sqrt{-\Delta}\}_{t>0}\\
   &\text{are $R$-bounded on}\quad L^{p}(\R^{n};X)\quad \text{for all}\quad p \in (p_{-},p_{+}).
\end{split}\right.\end{equation}
\begin{equation}\label{eq:calcAndKato}\left\{\begin{split}
\text{$L$ has an $H^\infty$ functional calculus in }L^{p}(\R^{n};X)\text{ and}\\
  \fNorm{\sqrt{L}u}{p}\eqsim\Norm{\nabla u}{p} \text{ for all } p \in (p_{-},p_{+}).
\end{split}\right.\end{equation}
\end{corollary}

\begin{remark}
It is interesting to note the connection between the above condition \eqref{eq:fourRbdd} and (a variant of) L.~Weis' characterization of so-called \emph{maximal regularity} \cite{Weis}: the \(R\)-boundedness of the set \(\{(I+t^2L)^{-1}\}_{t>0}\) in \(\bddlin(L^p(\R^n;X))\) is equivalent to the existence of a unique solution in \(L^q(\R;\domain(L))\cap W^{2,q}(\R;L^p(\R^n;X))\) of the problem \(-u''+Lu=f\) for each \(f \in L^q(\R;L^p(\R^n;X))\), where \(1<q<\infty\).
\end{remark}

We start with the proof of the corollary.

\begin{proof}
Let us first remark that the functional calculus in \eqref{eq:calcAndKato} implies the $R$-boundedness of \(\{(I+t^2L)^{-1}\}_{t>0}\) and
\(\{t\sqrt{L}(I+t^2L)^{-1}\}_{t>0}\) by \cite{kaltonweis} Theorem 5.3.
Using also the Kato estimates, we have that \(\eqref{eq:calcAndKato} \Rightarrow \eqref{eq:fourRbdd}\). Now consider a perturbed Hodge-Dirac operator \(\Pi_{B}\) with
$A_1=I$, $A_2=A$. Its resolvent can be computed as
\begin{equation*}
  (I-it\Pi_{B})^{-1}
  =\begin{pmatrix}
    (I+t^{2}L)^{-1} & -it(I+t^{2}L)^{-1}\diverg A\\
    it\nabla (I+t^{2}L)^{-1} & I+t^{2}\nabla(I+t^{2}L)^{-1}\diverg A
  \end{pmatrix}.
\end{equation*}

By Theorem \ref{basicmult}, $\nabla/\sqrt{-\Delta}$  is bounded from $L^p(\R^n;X)$ to $L^p(\R^n;X^n)$, and $\diverg/\sqrt{-\Delta}$
is bounded from $L^p(\R^n;X^n)$ to $L^p(\R^n;X)$. Using the boundedness of \(A\) on $L^p(\R^n;X^n)$,
the R-bisectoriality of $\Pi_B$ thus follows from \eqref{eq:fourRbdd}.
By Theorem~\ref{mainthm} the operator \(\Pi_{B}\) hence has an \(H^{\infty}\) functional calculus.
The functional calculus of \(L\) follows from the functional calculus of \(\Pi_{B}\) applied to functions of \(\Pi_{B}^{2}\).
The Kato estimates follow from  the functional calculus of \(\Pi_{B}\) applied to the sign function $z\mapsto z/\sqrt{z^2}$,
as in \cite{AKM} Corollary 2.11.
\end{proof}

Theorem \ref{mainthm} is a consequence of the following square function estimate, which is a vector-valued analogue of 
Proposition~4.8 in~\cite{AKM}. 

\begin{proposition}\label{mainprop}
Let \(X\) be a UMD Banach space such that both $X$ and $X^*$ have RMF. Consider perturbed Hodge-Dirac operators \(\Pi_{B}\) and \(\Pi_{B^{*}}\) in \(L^p(\R^n;X^N)\) for \(p\) in an open interval \((p_{-},p_{+})\subseteq (1,\infty)\).
Assume that \(\Pi_{B}\) and \(\Pi_{B^{*}}\) are \(R\)-bisectorial in \(L^p(\R^n;X^N)\) for all \(p \in (p_{-},p_{+})\).
Then we have
\begin{equation}\label{mainest}
   \sup_{1 \leq \abs{t} \leq 2}\Exp\BNorm{\sum_{k \in \Z} \radem_{k}Q_{2^{k}t}^{B}u}{
     L^p(\R^n;X^N)} \lesssim \Norm{u}{L^p(\R^n;X^N)},\quad
   \forall u \in \range(\Gamma).
\end{equation}
Moreover, the same estimates holds in \(L^p(\R^n;X^N)\) if the  triple \(\{\Gamma, B_{1}, B_{2}\}\) is replaced by \(\{\Gamma^{*}, B_{2}, B_{1}\}\),
and in \(L^{p'}(\R^n;(X^*)^N)\) if it is replaced by \(\{\Gamma, B_{1}^{*}, B_{2}^{*}\}\) or \(\{\Gamma^{*}, B_{2}^{*}, B_{1}^{*}\}\). 
\end{proposition}

This is proven in the rest of the paper. In fact, it suffices to prove the assertion with the triple \(\{\Gamma,B_1,B_2\}\), as written out in \eqref{mainest}, since the assumptions remain invariant when replacing this triple by any one of the three other possibilities. To simplify notation we will, moreover, only consider \(Q_{2^{k}}^{B}\) instead of \(Q_{2^{k}t}^{B}\), since the proofs remain the same in this generality.

We start our journey towards the proof of the main estimate \eqref{mainest} in the next section; in the rest of this section we show how to deduce Theorem \ref{mainthm} from Proposition \ref{mainprop}. We begin with the following:

\begin{lemma}[Hodge decomposition]\label{lem:Hodge}
Let $X$ be a reflexive Banach space, $1<p<\infty$, and $\Pi_B$ be a perturbed Hodge--Dirac operator which is bisectorial in $L^p(\R^n;X^N)$. Then the space decomposes as the following topological direct sum:
 \begin{equation*}
  L^{p}(\R^{n};X^N)
 = \kernel(\Pi_{B})\oplus \overline{\range(\Gamma)} \oplus
   \overline{\range(\Gamma^{*}_B)}.
\end{equation*}
\end{lemma}

\begin{proof}
On the abstract level, i.e., without making use of the structure of the Hodge--Dirac operators, the assumptions that $X$ (and then also $L^p(\R^n;X^N)$) is reflexive and $\Pi_B$ is bisectorial imply the decomposition
\begin{equation*}
  L^p(\R^n;X^N) = \kernel(\Pi_{B}) \oplus \overline{\range(\Pi_{B})}.
\end{equation*}
Moreover, the projection on
\(\overline{\range(\Pi_{B})}\) is given by
\begin{equation*}
  Pu = \lim_{t \to \infty} t^{2}\Pi_{B}^{2}(I+t^{2}\Pi_{B}^{2})^{-1}u.
\end{equation*}

In our specific situation, we further have the explicit formula
\begin{equation*}\begin{split}
& t^{2}\Pi_{B}^{2}(I+t^{2}\Pi_{B}^{2})^{-1} = \\
&  \begin{pmatrix}t^{2}A_{1}D^*A_{2}D(I+t^{2}A_{1}D^*A_{2}D)^{-1} & 0 \\ 0 & t^{2}D A_{1}D^* A_{2}(I+t^{2}D A_{1}D^* A_{2})^{-1}
  \end{pmatrix}.
\end{split}\end{equation*}
The projection \(P\) thus splits as \(P_{1}+P_{2}\), where \(P_i\) acts invariantly on \(L^p(\R^n;X^{n_i})\) and annihilates \(L^p(\R^n;X^{n_j})\) jor \(j\neq i\).
Since \(\overline{\range(\Gamma)} \subseteq\range(P) \cap L^{p}(\R^{n};X^{n_2}) = \range(P_2) \subseteq \overline{\range(\Gamma)}\), and
\(\overline{\range(\Gamma^{*}_{B})} \subseteq \range(P) \cap L^{p}(\R^{n_1};X) = \range(P_1) \subseteq \overline{\range(\Gamma^{*}_{B})}\),
this gives the Hodge decomposition.
\end{proof}

\begin{proof}[of Theorem \ref{mainthm}]
The fact that \eqref{eq:PiBfc} \(\Rightarrow\) \eqref{eq:PiBRsect} is essentially contained in \cite{kaltonweis}, Theorem 5.3, where it is stated for sectorial (rather than bisectorial) operators.
Likewise, the equivalence between the square function estimates
\begin{equation}\label{quadest}
\left\{\begin{split}
  &\sup_{1 \leq\abs{t}\leq 2}
   \mathbb{E}\BNorm{\sum_{k \in \Z} \varepsilon_{k}Q_{2^{k}t}^{B}u}{p} \lesssim \Norm{u}{p} \quad
   \forall u \in L^{p}(\R^{n};X^N), \\
  &\sup_{1 \leq\abs{t}\leq 2}
   \mathbb{E}\BNorm{\sum_{k \in \Z} \varepsilon_{k}(Q_{2^{k}t}^{B})^{*}u}{p'} \lesssim \Norm{u}{p'} \quad
   \forall u \in L^{p'}(\R^{n};(X^*)^{N}),
\end{split}\right.
\end{equation}
and the functional calculus of \(\Pi_{B}\) is proven in \cite{kunstweis} Theorem 12.17 for sectorial operators but the proof carries over to the bisectorial situation.

We thus have to show that \eqref{eq:PiBRsect} implies \eqref{quadest}. By Lemma~\ref{lem:Hodge}, it suffices to do this separately for $u$ in each of the three components of the Hodge decomposition.  Now \(Q_{2^{k}}^{B}u = 0\) for all \(u \in \kernel(\Pi_{B})\), and Proposition \ref{mainprop} gives the first estimate in (\ref{quadest}) for \(u \in \range(\Gamma)\).
On \(\range(\Gamma^{*}_B)\), we then apply Proposition \ref{mainprop} with the triple \((\Gamma,B_{1},B_{2})\) replaced by
\((\Gamma^{*},B_{2},B_{1})\).

This gives
\begin{equation*}
  \sup_{1 \leq \abs{t} \leq 2} \mathbb{E}\BNorm{\sum_{k \in \Z}
   \varepsilon_{k} 2^{k}tB_{2}\Gamma B_{1}(I+(2^{k}t
    (\Gamma^{*}+B_{2}\Gamma B_{1}))^{2})^{-1}u}{p} \lesssim \Norm{u}{p},
\end{equation*}
for all  \(u \in R(\Gamma^{*})\); by simple manipulation, this is equivalent to
\begin{equation*}
  \sup_{1 \leq \abs{t} \leq 2} \mathbb{E}\BNorm{\sum \limits _{k \in \Z}
    \varepsilon_{k}2^{k}t\Gamma (I+(2^{k}t\Pi_{B})^{2})^{-1}B_{1}u}{p} \lesssim \Norm{u}{p}
  \quad \forall u \in R(\Gamma^{*}),
\end{equation*}
and then in turn to
\begin{equation*}
  \sup_{1 \leq \abs{t} \leq 2} \mathbb{E}\BNorm{
     \sum_{k \in \Z} \varepsilon_{k}Q_{2^{k}t}^{B}u}{p} \lesssim \Norm{u}{p} \quad
   \forall u \in \range(\Gamma^{*}_B),
\end{equation*}
since $\range(\Gamma_B^*)=B_1\range(\Gamma^*)$. 

To obtain the dual estimates, one remarks that the above reasoning can be applied to \(\Pi_{B}^{*} = \Gamma^{*} + B_{2}^{*}\Gamma B_{1}^{*}\) and
\(\Pi_{B^{*}}^{*} = \Gamma + B_{1}^{*}\Gamma^{*} B_{2}^{*}\). Indeed, these operators are $R$-bisectorial on \(L^{p'}(\R^{n};(X^{*})^N)\) by the duality of $R$-bounds (Lemma 3.1 in \cite{kaltonweis}; here one needs the fact that UMD spaces have nontrivial type.)
\end{proof}

\begin{remark}
The reader familiar with Hodge--Dirac operators will have noticed the special form of our operators \(\Gamma\) and 
\(\Gamma^{*}\), and, in particular, the fact that we are not working at the level of generality of \cite{AKM}.
However, the proof of Proposition~\ref{mainprop} carries over to the following situation.

\begin{proposition}
Let \(X\) be a UMD Banach space such that both $X$ and $X^*$ have RMF. 
Let \(\Gamma\) be a nilpotent first order differential operator with constant coefficients in \(\mathcal{L}(\C^{N})\) satisfying \(\Pi^{3} = -\Delta \Pi\),
where \(\Pi = \Gamma+\Gamma^{*}\).
Let \(B_{1},B_{2} \in L^{\infty}(\R^{n};\mathcal{L}(\C^{N}))\) be such that \(\Gamma^{*}B_{2}B_{1}\Gamma^{*} = 0 = \Gamma B_{1}B_{2}\Gamma\).
Assume that \(\Pi_{B} = \Gamma +B_{1}\Gamma^{*}B_{2}\) is R-bisectorial on \(L^{q}(\R^{n};X^{N})\) for all \(q\in(p-\varepsilon,p+\varepsilon)\), where \(\varepsilon >0\).
Then we have
\begin{equation*}\begin{split}
   \sup_{1 \leq \abs{t} \leq 2}\Exp\BNorm{\sum_{k \in \Z} \radem_{k}2^{k}t\Pi_{B}(I+(2^{k}t\Pi_{B})^{2})^{-1}u}{
     L^p(\R^n;X^N)} \lesssim &\Norm{u}{L^p(\R^n;X^N)}, \\
   &\forall u \in \range(\Gamma).
\end{split}\end{equation*}
\end{proposition}

This holds, in particular, in the case where \(\Gamma\) is an exterior derivative.
However, the \(L^{p}\) Hodge decomposition of Lemma \ref{lem:Hodge} is no more automatic in this situation.
To deduce a version of Theorem \ref{mainthm} in this more general setting one would thus need to have the existence of the Hodge decomposition as an assumption. 
Since our main focus is the original square root problem, we chose not to work in this generality in order to keep the paper more readable.  The \(L^{p}\) theory of more general Hodge--Dirac operators will be considered elsewhere.
\end{remark}


\section{Miscellaneous propositions}\label{sec:props}

This section is a sm\"org{\aa}sbord of vector-valued analogues of a number of classical estimates of Analysis, which we need in the subsequent developments. We start with a vector-valued version of the Poincar\'e inequality. Below, \(u\cdot v\) denotes the dot product of \(u,v \in \R^{n}\), \(\tau_{h}\) stands for the translation operator defined by \(\tau_{h}f(x)=f(x+h)\), and \(1_{Q}\) denotes the characteristic function of the set \(Q\).

\begin{proposition}[Poincar\'e inequality]\label{poincare}
Let \(X\) be a Banach space, and \(1\leq p<\infty\). For \(u \in W^{1,p}(\R^{n};X)\),  and \(m \in \Z^{n}\)  we have
\begin{equation*}\begin{split}
  \Exp &\BNorm{\sum_{k\in\Z}\radem_k\sum_{Q\in\triangle_{2^k}}
    1_Q(u_k-\ave{u_k}_{Q+2^k m})}{L^p(\R^n;X)} \\
  &\lesssim\int_{[-1,1]^n}\int_0^1
  \Exp\BNorm{\sum_{k\in\Z}\radem_k 2^k(m+z)\cdot\nabla\tau_{t2^k(m+z)}u_k
   }{L^p(\R^n;X)}\ud t \ud z.
\end{split}\end{equation*}
\end{proposition}

\begin{proof}
For \(x\in Q\in\triangle_{2^k}\), we observe that \(Q\subset x+2^k[-1,1]^n\). Hence
\begin{equation*}\begin{split}
  &u_k(x)-\ave{u_k}_{Q+2^k m}\\
  &=\int_{[-1,1]^n}\big[u_k(x)-u_k(x+2^k(m+z))\big]1_Q(x+2^k z)\ud z \\
  &=\int_{[-1,1]^n}\int_0^1 -2^k(m+z)\cdot\nabla u_k(x+t2^k(m+z))\ud t\,
    1_Q(x+2^k z)\ud z.
\end{split}\end{equation*}
The assertion follows after bringing the integrals outside the norm and discarding the indicators \(1_Q(x+2^k z)\) by the contraction principle 
 \ref{contraction-principle} .
\end{proof}

Here is a useful Banach space version of another classical inequality:

\begin{proposition}[Schur's estimate]\label{prop:Schur}
Let $\mathcal{X}$, $\mathcal{Y}$ and $\mathcal{Z}$ be Banach spaces, the last two with finte cotype. For $i,j\in\Z$, let $\alpha(i,j)$ be positive numbers satisfying
\begin{equation*}
  \sup_i\sum_j\alpha(i,j)\lesssim 1,\qquad
  \sup_j\sum_i\alpha(i,j)\lesssim 1,
\end{equation*}
and let $T_{i,j}\in\bddlin(\mathcal{Y},\mathcal{Z})$, $D_i\in\bddlin(\mathcal{X},\mathcal{Y})$ be operators satisfying
\begin{equation*}
  \rbound(\frac{1}{\alpha(i,j)}T_{i,j}:
   i,j\in\Z)\lesssim 1,\qquad
  \Exp\BNorm{\sum_i\radem_i D_i x}{\mathcal{Y}}\lesssim\Norm{x}{\mathcal{X}}
\end{equation*}
for all $x\in \mathcal{X}$. Then there holds
\begin{equation*}
  \Exp\BNorm{\sum_{i,j}\radem_j T_{i,j}D_i x}{\mathcal{Z}}
  \lesssim\Norm{x}{\mathcal{X}}.
\end{equation*}
\end{proposition}

\begin{proof}
Under the assumption that $\mathcal{Y}$ and $\mathcal{Z}$ have finite cotype, we may replace
the Rademacher-variables $\radem_j$ in the assumptions and the
claim by independent standard Gaussian random variables
$\gamma_j$ by \eqref{eq:RadVsGauss}. We write the left-hand side of the modified assertion
as
\begin{equation*}
  \Exp\BNorm{\sum_j\gamma_j\sum_i\alpha(i,j)^{1/2}
   \frac{1}{\alpha(i,j)}T_{i,j}\alpha(i,j)^{1/2}D_i x}{\mathcal{Z}}.
\end{equation*}
Then, as in \cite{HytPot} Proposition 2.1, let
\begin{equation*}
  x_{i,j}:=\frac{1}{\alpha(i,j)}T_{i,j}\alpha(i,j)^{1/2}D_i x,\qquad
  y_{j} := \sum \limits _{i} \alpha(i,j)^{\frac{1}{2}}x_{i,j}.
\end{equation*}
For \(x^{*}\in\mathcal{X}^{*}\), we have 
\begin{equation*}
  \sum_{j} \abs{\pair{y_{j}}{x^{*}}}^{2}
  \leq \sup_j\Big(\sum_{i} \alpha(i,j)\Big)
  \sum_{i,j} \abs{\pair{x_{i,j}}{x^{*}}}^{2}.
\end{equation*}
Now Proposition 3.7 in \cite{pisier:facto} states that
\begin{equation*}\begin{split}
  \sum_{j} \abs{\pair{y_{j}}{x^{*}}}^{2} &\leq C^2   \sum_{i,j} \abs{\pair{x_{i,j}}{x^{*}}}^{2}\quad\forall x^{*} \in X^{*}\\
  \Rightarrow\quad \Exp \|\sum_{j} \gamma_{j}y_{j}\| &\leq C \Exp \|\sum_{i,j} \gamma_{i,j}x_{i,j}\|,
\end{split}\end{equation*}
where \((\gamma_{i,j})_{i,j\in\Z}\) is a double-indexed sequence
of independent standard Gaussian variables.
Therefore, using our \(R\)-boundedness assumption, we have
\begin{equation*}\begin{split}
 &\Exp\BNorm{\sum_j\gamma_j\sum_i\alpha(i,j)^{1/2}
   \frac{1}{\alpha(i,j)}T_{i,j}\alpha(i,j)^{1/2}D_i x}{\mathcal{Z}}\\
&  \leq\sup_j\Big(\sum_i\alpha(i,j)\Big)^{1/2}
  \Exp\BNorm{\sum_{i,j}\gamma_{i,j}\frac{1}{\alpha(i,j)}T_{i,j}
    \alpha(i,j)^{1/2}D_i x}{\mathcal{Z}}\\
&  \lesssim\Exp\BNorm{\sum_{i,j}\gamma_{i,j}
    \alpha(i,j)^{1/2}D_i x}{\mathcal{Y}}.
\end{split}\end{equation*}
By reorganization, the last expression is equal to
\begin{equation*}
  \Exp\BNorm{\sum_i\Big(\sum_j\alpha(i,j)^{1/2}\gamma_{i,j}\Big)
    D_i x}{\mathcal{Y}}
  =:\Exp\BNorm{\sum_i\tilde\gamma_i D_i x}{\mathcal{Y}}.
\end{equation*}
By basic properties of Gaussian sums, the random variables
\(\tilde\gamma_i\) are again independent Gaussian, with variance
\begin{equation*}
  \Exp\tilde\gamma_i^2=\sum_j\alpha(i,j)\lesssim 1.
\end{equation*}
By the contraction principle  \ref{contraction-principle} , the random sum with
\(\tilde\gamma_i\)'s is then dominated by a random sum with standard
Gaussian variables, and using the assumption on the operators
\(D_i\) we complete the argument.
\end{proof}

In the rest of this section, we make use of the Haar system of functions. Recall that in $\R^n$ there are $2^n-1$ Haar functions $h_Q^{\eta}$, $\eta\in\{0,1\}^n\setminus\{0\}$, associated with every dyadic cube $Q\in\triangle$. For our purposes, it is most convenient to normalize them in $L^{\infty}(\R^n)$ so that $\abs{h_Q^{\eta}}=1_Q$. We often need only one (say, the ``first'') of the $h_Q^{\eta}$ for each $Q$, and so we adopt the notation  $h_Q:=h_{Q}^{(1,0,...,0)}:=1_{Q_+}-1_{Q_-}$, where $Q_+$ and $Q_-$ are two halves of $Q$.

\begin{lemma}[Sign-invariance]\label{lem:signInvar}
Let $X$ be any Banach space, $1\leq p<\infty$ and $u_Q\in L^p(\R^n;X)$ for all $Q\in\triangle$. Then 
\begin{equation*}\begin{split}
  \Exp\BNorm{\sum_k\radem_k\sum_{Q\in\triangle_{2^k}}
    1_Q u_Q}{L^{p}(\R^{n};X)}
  &\eqsim\Exp\BNorm{\sum_k\radem_k\sum_{Q\in\triangle_{2^k}}
    h_Q u_Q}{L^{p}(\R^{n};X)} \\
  &\eqsim\Exp\BNorm{\sum_{Q\in\triangle}
    \radem_Q 1_Q u_Q}{L^{p}(\R^{n};X)}.
\end{split}\end{equation*}
\end{lemma}

\begin{proof}
Using Kahane's inequality \ref{kk-inequalities}, and Fubini's Theorem, we have
\begin{equation*}\begin{split}
  &\Exp\BNorm{\sum_k\radem_k\sum_{Q\in\triangle_{2^k}}
    1_Q u_Q}{L^{p}(\R^{n};X)} \\
  &\eqsim \Big( \int \limits _{\R^{n}}   \Exp\BNorm{\sum_k\radem_k\sum_{Q\in\triangle_{2^k}}
    1_Q(y) u_Q(y)}{X}^{p}dy \Big)^{1/p}
\end{split}\end{equation*}
For a fixed \(y \in \R^{n}\), and a scale \(k \in \Z\), there exists a unique dyadic cube \(Q_{k,y} \in \triangle_{2^{k}}\) containing \(y\).
Therefore, by the contraction principle \ref{contraction-principle}
\begin{equation*}
 \Exp\BNorm{\sum_k\radem_k\sum_{Q\in\triangle_{2^k}}
    1_Q(y) u_Q(y)}{X} \simeq  \Exp\BNorm{\sum_k\radem_k\sum_{Q\in\triangle_{2^k}}
    h_Q(y) u_Q(y)}{X}.
\end{equation*}
This gives the first equivalence. A similar argument applies to the second.
\end{proof}

We next recall a result of Figiel from \cite{Figiel:88}. Our need for it is no surprise, since it is also a fundamental ingredient in Figiel's vector-valued $T(1)$ theorem \cite{Figiel:90}.

\begin{proposition}[Figiel]\label{prop:Figiel}
Let \(X\) be a UMD Banach space, and $1<p<\infty$. Then for all \(m \in \Z^{n}\) and $x_Q^{\eta}\in X$
\[
  \BNorm{\sum_{Q\in\triangle}\sum_{\eta}
    x_Q^{\eta} h_{Q+\ell(Q)m}^{\eta}}{L^p(\R^n;X)}
  \lesssim \log(2+\abs{m})\BNorm{\sum_{Q\in\triangle}\sum_{\eta}
    x_Q^{\eta} h_Q^{\eta}}{L^p(\R^n;X)}.
\]
\end{proposition}

\begin{corollary}\label{cor:Figiel}
Let \(X\) be a UMD Banach space, and \(1<p<\infty\). For \((u_k)_{k \in \Z} \subset L^{p}(\R^{n};X)\), we have
\[
  \Exp\BNorm{\sum_k\radem_k\sum_{Q\in\triangle_{2^k}}
    1_{Q+2^k m}\ave{u_k}_Q}{L^p(\R^n;X)}
  \lesssim \log(2+\abs{m})\BNorm{\sum_k\radem_k u_k}{
    L^p(\R^n;X)}
\]
\end{corollary}

\begin{proof}
By sign-invariance,  and unconditionality of the Haar system, $\log(2+\abs{m})^{-1}$ times the left-hand side is equivalent to
\begin{equation*}\begin{split}
  \frac{1}{\log(2+\abs{m})}\Exp &\BNorm{\sum_k\radem_k\sum_{Q\in\triangle_{2^k}}
    h_{Q+2^k m}\ave{u_k}_Q}{}
  \lesssim \Exp\BNorm{\sum_k\radem_k
    \sum_{Q\in\triangle_{2^k}}1_Q\ave{u_k}_Q}{}\\
  &= \Exp\BNorm{\sum_k\radem_k A_{2^k}u_k}{}
  \lesssim\Exp\BNorm{\sum_k u_k}{},
\end{split}\end{equation*}
where Stein's inequality \ref{rbdd-At} was used in the last step.
\end{proof}

The following Lemma, too, is closely related to Proposition~\ref{prop:Figiel}, but unlike in the easy Corollary above, we now have to employ the techniques of Figiel's proof \cite{Figiel:88} rather than just his result. Similar martingale arguments inspired by \cite{Figiel:88} were also recently used in~\cite{Hytonen:Tb}.

\begin{lemma}\label{lem:martingale}
Let $X$ be a UMD space, and $1<p<\infty$. Let further $k\in\Z_+$, $\ell\in\{0,\ldots,k\}$, and $x_Q\in X$ for all $Q\in\triangle$. For each $Q\in\triangle$, let  $E(Q),F(Q)\subset Q$ be two disjoint subsets such that: both $E(Q)$ and $F(Q)$ are unions of some dyadic cubes $R\in\triangle_{2^{-k}\ell(Q)}$, and $\abs{F(Q)}\leq\abs{E(Q)}$.
Then
\begin{equation*}\begin{split}
  &\Big(\Exp\BNorm{\sum_{j\equiv\ell}\radem_j
    \sum_{Q\in\triangle_{2^j}}1_{F(Q)} x_Q}{L^p(\R^n;X)}^p\Big)^{1/p} \\
  &\lesssim_{p,X}\Big(\Exp\BNorm{\sum_{j\equiv\ell}\radem_j
    \sum_{Q\in\triangle_{2^j}}1_{E(Q)} x_Q}{L^p(\R^n;X)}^p\Big)^{1/p},
\end{split}\end{equation*}
where $j\equiv\ell$ is shorthand for $j\equiv\ell\mod (k+1)$.
\end{lemma}

\begin{proof}
Let
\begin{equation*}
  E(Q)=\bigcup_{i=1}^{I(Q)}R_i(Q),\qquad
  F(Q)=\bigcup_{i=1}^{J(Q)}S_i(Q),
\end{equation*}
where $R_i(Q),S_i(Q)\in\triangle_{2^{-k}\ell(Q)}$, the unions are disjoint, and therefore $J(Q)\leq I(Q)\leq 2^{kn}$ by assumption. Writing 
$1_{F(Q)} = \sum_{i} 1_{S_i(Q)}$, $1_{E(Q)} = \sum_{i} 1_{R_i(Q)}$, and using sign-invariance, the claim is seen to be equivalent to
\begin{equation}\label{eq:martToProve}\begin{split}
  &\Exp\BNorm{\sum_{j\equiv\ell}\radem_j
    \sum_{Q\in\triangle_{2^j}}\sum_{i=1}^{J(Q)}h_{S_i(Q)} x_Q}{L^p(\R^n;X)}^p \\
  &\lesssim\Exp\BNorm{\sum_{j\equiv\ell}\radem_j
    \sum_{Q\in\triangle_{2^j}}\sum_{i=1}^{I(Q)}h_{R_i(Q)} x_Q}{L^p(\R^n;X)}^p.
\end{split}\end{equation}

We may consider the point in our probability space being fixed for a while, so that the $\radem_j$ are just some given signs. For each $j\equiv\ell$ and $Q\in\triangle_{2^j}$, we introduce auxiliary functions as follows:
\begin{alignat*}{2}
  &d^{\pm 1}_{Q,i}:=\radem_j\frac{1}{2}(h_{R_i(Q)}\pm h_{S_i(Q)})x_Q,\qquad & 1\leq i\leq J(Q),\\
  &d^0_{Q,i}:=\radem_jh_{R_i(Q)}x_Q,\qquad & J(Q)<i\leq I(Q),
\end{alignat*}
and finally
\begin{equation*}
   d^{\pm 1}_j :=\sum_{Q\in\triangle_{2^j}}\sum_{i=1}^{J(Q)}d^{\pm 1}_{Q,i},\qquad
   d^{0}_j :=\sum_{Q\in\triangle_{2^j}}\sum_{i=J(Q)+1}^{I(Q)}d^{0}_{Q,i}.
\end{equation*}

Let us make a key observation. If $Q,Q'\in\triangle$ appear in the claimed estimate \eqref{eq:martToProve} and $\ell(Q)>\ell(Q')$, then $\ell(Q)\geq 2^{k+1}\ell(Q')$. The functions $d^{\theta}_{Q,i}$ are constant on halves of dyadic cubes of side-length $2^{-k}\ell(Q)$, and hence they are constants on $Q'$.

We now define the following $\sigma$-algebras:
\begin{equation*}\begin{split}
  \mathscr{F}_j^0 &:=\sigma(\triangle_{2^{j-1}}),\\
  \mathscr{F}_j^1 &:=\sigma(\mathscr{F}_j^0,\{d^{+1}_{Q,i}:Q\in\triangle_{2^j},1\leq i\leq J(Q)\}),\\
  \mathscr{F}_j^2 &:=\sigma(\mathscr{F}_j^1,\{d^{-1}_{Q,i}:Q\in\triangle_{2^j},1\leq i\leq J(Q)\}),
\end{split}\end{equation*}
where \(\sigma(S)\) denotes the sigma algebra generated by the elements of \(S\), and 
\(\{d^{\pm1}_{Q,i}:Q\in\triangle_{2^j},1\leq i\leq J(Q)\}\) denotes the sets, indexed by \(Q\in\triangle_{2^j}\) and \(i\), of sets \((d^{\pm1}_{Q,i})^{-1}(B)\) where \(B \subset \R\) is a Borelian set.
Then
\begin{equation*}
  \ldots\subseteq\mathscr{F}_{\ell+\nu(k+1)}^0
  \subseteq\mathscr{F}_{\ell+\nu(k+1)}^1
  \subseteq\mathscr{F}_{\ell+\nu(k+1)}^2
  \subseteq\mathscr{F}_{\ell+(\nu-1)(k+1)}^0\subseteq\ldots
\end{equation*}
is a filtration of $\R^n$ which generates the Borel $\sigma$-algebra, and
\begin{equation*}
  \ldots,\quad d^0_{\ell+(\nu+1)(k+1)},\quad d^{+1}_{\ell+\nu(k+1)},\quad
  d^{-1}_{\ell+\nu(k+1)},\quad d^0_{\ell+\nu(k+1)},\quad\ldots
\end{equation*}
is a martingale difference sequence, with respect to this filtration.

By the very definition of UMD spaces, there holds
\begin{equation*}
  \BNorm{\sum_{j\equiv\ell}\sum_{\theta\in\{0,\pm 1\}}
    \theta d^{\theta}_j}{L^p(\R^n;X)}
  \lesssim \BNorm{\sum_{j\equiv\ell}\sum_{\theta\in\{0,\pm 1\}}
    d^{\theta}_j}{L^p(\R^n;X)}.
\end{equation*}
But it is immediate to see that this estimate, after taking the expectation with respect to the $\radem_j$ on both sides, is precisely the desired inequality \eqref{eq:martToProve}.
\end{proof}

\begin{remark}
In the above lemma, the disjointness assumption for \(E(Q)\) and \(F(Q)\) can be dropped. 
Writing \(1_{F(Q)}\) as \(1_{F(Q)}.1_{E(Q)}+1_{F(Q)\backslash E(Q)}\), one can apply the above proof with \(F(Q)\) replaced by \(F(Q)\backslash E(Q)\), and handle the other term using sign-invariance and the contraction principle.
\end{remark}


\section{Vector-valued inequalities for the unperturbed operator}\label{unpert}

For the unperturbed operator $\Pi$, we define $R_t$, $P_t$ and $Q_t$ by simply dropping the $B$'s from the formulae~\eqref{eq:resolvents}. We also set
\begin{equation*}
  \mathscr{P}_t=(I-t^2\Delta)^{-1},\qquad
  \mathscr{Q}_t=t\nabla\mathscr{P}_t,\qquad
  \mathscr{Q}_t^*=-t\mathscr{P}_t\diverg;
\end{equation*}
as it turns out, the assumption \eqref{eq:DDstar} often helps to reduce the more complicated Hodge-Dirac resolvents to this canonical family of operators. Note that
\begin{equation*}
  \mathscr{Q}_t^*\mathscr{Q}_t=-t^2\Delta(I-t^2\Delta)^{-2}.
\end{equation*}

An important component of our work is the analogue  between the harmonic and the dyadic worlds, and in particular the idea that $\mathscr{P}_t$ and $A_t$ are roughly the same. This heuristic will be quantified and proved later on.

\begin{proposition}
\label{funcalc}
Let \(X\) be a UMD Banach space and \(1<p<\infty\). Then the Hodge-Dirac operator \(\Pi\) has an $H^{\infty}(S_{\theta})$ functional calculus
on  \(L^{p}(\R^{n};X^N)\) for every $\theta>0$.
\end{proposition}

\begin{proof}
With the help of the Fourier transform and the elementary functional calculus of selfadjoint matrices, the functional calculus of $\Pi$ may be computed explicitly. In fact, it follows from the assumption \eqref{eq:DDstar} that the symbol $\hat{\Pi}(\xi)$ of the differential operator $\Pi$ satisfies $\hat{\Pi}(\xi)^3=\abs{\xi}^2\hat{\Pi}(\xi)$, which implies that the only possible eigenvalues of the matrix $\hat{\Pi}(\xi)$ are $0$ and $\pm\abs{\xi}$. Functions of such matrices are readily computed, and transforming back we find that
\begin{equation}\label{eq:fcPi}
  f(\Pi)=f_o(\sqrt{-\Delta})\frac{\Pi}{\sqrt{-\Delta}}
  +[f_e(\sqrt{-\Delta})-f(0)]\frac{\Pi^2}{-\Delta}
  +f(0)I,
\end{equation}
where $f_o(z):=\frac{1}{2}\big(f(z)-f(-z)\big)$ and $f_e(z):=\frac{1}{2} \big(f(z)+f(-z)\big)$ are the odd and even parts of $f$, respectively. All the operators above are Fourier multipliers, whose boundedness on $L^p(\R^n;X^N)$ follows from the Multiplier Theorem~\ref{basicmult}.
\end{proof}

Note that \eqref{eq:fcPi} and $\Pi^3=-\Delta\Pi$ imply in particular that
\begin{equation}\label{eq:fcEven}
  g(\Pi^2)\Pi=g(-\Delta)\Pi,
\end{equation}
i.e., on $\range(\Pi)$ the functional calculus of $\Pi^2$ is just the functional calculus of $-\Delta$.

\begin{lemma}\label{lem:multiplier}
Let $X$ be a UMD space, $1<p<\infty$, and $2M>n+1$. For $z\in\R^n$, $t\in[0,1]$, and $u\in L^p(\R^n;X)$, there holds
\begin{equation*}
  \Exp\BNorm{\sum_k\radem_k 2^k z\cdot
    \nabla\tau_{t2^k z}\mathscr{P}_{2^k}^M u}{L^p(\R^n;X)}
  \lesssim(1+\abs{z})^{n+1}\Norm{u}{L^p(\R^n;X)}.
\end{equation*}
\end{lemma}

\begin{proof}
The function inside the norm on the left is a Fourier multiplier transformation of $u$ with the symbol
\begin{equation*}
  \sigma(\xi)=\sum_k\radem_k 2^k z\cdot i\xi\, e^{it2^kz\cdot\xi}\cdot
  (1+2^{2k}\abs{\xi}^2)^{-M}.
\end{equation*}
For every $\alpha\in\{0,1\}^n$, a straightforward computation shows, given the assumption $2M>n+1$, that
\begin{equation*}
  \abs{\xi}^{\abs{\alpha}}\abs{D^{\alpha}\sigma(\xi)}
  \lesssim (1+\abs{z})^{1+\abs{\alpha}}\lesssim (1+\abs{z})^{1+n}.
\end{equation*}
The assertion hence follows from the Multiplier Theorem~\ref{basicmult}.
\end{proof}

\begin{lemma}\label{PtVsPtN}
Let \(X\) be a UMD Banach space, \(1<p<\infty\), and $M\in\Z_+$. For \(u \in L^{p}(\R^{n};X)\) we have
\begin{equation*}
   \Exp\BNorm{\sum_{k\in\Z}\radem_k
     (\mathscr{P}_{2^k}-\mathscr{P}_{2^k}^M)u}{L^p(\R^n;X^N)}
   \lesssim \Norm{u}{L^p(\R^n;X^N)}
\end{equation*}
\end{lemma}

\begin{proof}
This is a Fourier multiplier estimate again. One may either directly study the multiplier on the left like in Lemma~\ref{lem:multiplier}, or argue in a slightly more step-by-step fashion as follows: Observe first that $\mathscr{P}_t^{j-1}-\mathscr{P}_t^j=-t^2\Delta\mathscr{P}_t^j
=\mathscr{P}_t^{j-2} \mathscr{Q}_t^*\mathscr{Q}_t$ for all $j=2,\ldots,N$. The symbols onf $\mathscr{P}_t$ have uniformly bounded variation, so the operators are $R$-bounded by Proposition~\ref{RbddMult}, and thus
\begin{equation*}\begin{split}
  \Exp\BNorm{\sum_{k\in\Z}\radem_k
     (\mathscr{P}_{2^k}-\mathscr{P}_{2^k}^M)u}{}
  &\leq\sum_{j=2}^N\Exp\BNorm{\sum_{k\in\Z}\radem_k
     (\mathscr{P}_{2^k}^{j-1}-\mathscr{P}_{2^k}^j)u}{}\\
  &\lesssim\sum_{j=2}^N\Exp\BNorm{\sum_{k\in\Z}\radem_k
     \mathscr{Q}_{2^k}^*\mathscr{Q}_{2^k} u}{}
  \lesssim\Norm{u}{},
\end{split}\end{equation*}
where the final quadratic estimate again follows from the Multiplier Theorem~\ref{basicmult}.
\end{proof}

We have now accumulated enough knowledge to prove the following estimate showing that $\mathscr{P}_t$ is almost like its average $A_t\mathscr{P}_t$, in the precise sense of the quadratic estimate. In the rest of this section we are going to show the ``dual'' property that also $A_t$ is almost like $A_t\mathscr{P}_t$, thus justifying our heuristic of the ``equivalence'' of $A_t$ and $P_t$. 

\begin{proposition}\label{prop:APminusP}
Let $X$ be a UMD space, and $1<p<\infty$. Then for all $u\in L^p(\R^n;X)$, there holds
\begin{equation*}
  \Exp\BNorm{\sum_k\radem_k(A_{2^k}-I)\mathscr{P}_{2^k}u}{
    L^p(\R^n;X)}
  \lesssim\Norm{u}{L^p(\R^n;X)}.
\end{equation*}
\end{proposition}

\begin{proof}
Since the operators $A_{2^k}-I$ are $R$-bounded, and the differences $\mathscr{P}_{2^k}-\mathscr{P}_{2^k}^M$ satisfy the quadratic estimate of Lemma~\ref{PtVsPtN} (taking \(2M>n+1\)), it suffices to prove the claim with $\mathscr{P}_{2^k}$ replaced by $\mathscr{P}_{2^k}^M$. The left side of the modified claim is
\begin{equation*}\begin{split}
  &\Exp\BNorm{\sum_k\radem_k\sum_{Q\in\triangle_{2^k}}
    1_Q(\mathscr{P}_{2^k}^M u
   -\ave{\mathscr{P}_{2^k}^M u}_Q)}{L^p(\R^n;X)} \\
  &\lesssim\int_{[-1,1]^n}\int_0^1
    \Exp\BNorm{\sum_k\radem_k 2^k z\cdot\nabla
    \tau_{t2^k z}\mathscr{P}_{2^k}^M u}{L^p(\R^n;X)}\ud t\ud z \\
  &\lesssim\int_{[-1,1]^n}\int_0^1(1+\abs{z})^{n+1}
     \Norm{u}{L^p(\R^n;X)}\ud t\ud z
   \lesssim\Norm{u}{L^p(\R^n;X)}
\end{split}\end{equation*}
by Proposition~\ref{poincare} and Lemma~\ref{lem:multiplier}.
\end{proof}

 Our next Proposition is a vector-valued analogue of 
Proposition 5.7 in \cite{AKM}. 

\begin{proposition}\label{prop:APminusA}
Let \(X\) be a UMD space, and \(1<p<\infty\). For \(u \in L^{p}(\R^{n};X)\), we have
\begin{equation*}
  \Exp\BNorm{\sum_j\radem_jA_{2^j}(\mathscr{P}_{2^j}-I)u}{L^p(\R^n;X)}
  \lesssim\Norm{u}{L^p(\R^n;X)}.
\end{equation*}
\end{proposition}

\begin{proof}
As a preparation, observe that $\sum_{i\in\Z}\mathscr{Q}_{2^i}^*\mathscr{Q}_{2^i}$ is represented by the Fourier multiplier 
$\sum_{i\in\Z}(2^{i}\abs{\xi})^2(1+(2^i\abs{\xi})^2)^{-2}$ which, as well as its reciprocal, satisfies the conditions of the Multiplier Theorem~\ref{basicmult}. This implies the two-sided estimate
\begin{equation*}
  \BNorm{\sum_i \mathscr{Q}_{2^i}^*\mathscr{Q}_{2^i}
     u}{L^p(\R^n;X)}\eqsim\Norm{u}{L^p(\R^n;X)}.
\end{equation*}

Thus, it suffices to prove
\begin{equation}\label{eq:APAtoProve}
  \Exp\BNorm{\sum_{i,j}\radem_j A_{2^j}(\mathscr{P}_{2^j}-I)
    \mathscr{Q}_{2^i}^*
    \mathscr{Q}_{2^i}u}{L^p(\R^n;X)}\lesssim\Norm{u}{L^p(\R^n;X)}.
\end{equation}
Since also
\begin{equation*}
  \Exp\BNorm{\sum_i\radem_i\mathscr{Q}_{2^i}u}{L^p(\R^n;X)}
  \lesssim\Norm{u}{L^p(\R^n;X^n)},
\end{equation*}
again by Theorem~\ref{basicmult} (say), \eqref{eq:APAtoProve} will follow from Schur's estimate~\ref{prop:Schur}  (with $\mathcal{X}=\mathcal{Z}=L^p(\R^n;X)$, and 
$\mathcal{Y}=L^p(\R^n;X^n)$),  once we show that
\begin{equation}\label{Rbd-for-Schur}
  \rbound(2^{\delta\abs{j-i}}A_{2^j}(\mathscr{P}_{2^j}-I)
  \mathscr{Q}_{2^i}^*:i,j\in\Z)\lesssim 1
\end{equation}
for some $\delta>0$. Since  $(I-\mathscr{P}_t)\mathscr{Q}_s^*=\frac{t}{s}(I-\mathscr{P}_s)\mathscr{Q}_t^*$ and
$\mathscr{P}_t\mathscr{Q}_s^*=\frac{s}{t}\mathscr{P}_s\mathscr{Q}_t^*$  for all $s,t>0$, and all the families
$A_{2^j}$, $\mathscr{P}_{2^j}$ and $\mathscr{Q}_{2^j}^*$, $j\in\Z$, are $R$-bounded on the relevant spaces, it is immediate that
\begin{equation*}\begin{split}
  &\rbound(2^{i-j}A_{2^j}(\mathscr{P}_{2^j}-I)
     \mathscr{Q}_{2^i}^*:i\geq j)
  =\rbound(A_{2^j}\mathscr{Q}_{2^j}^*
      (\mathscr{P}_{2^i}-I):i\geq j)\lesssim 1,\\
  &\rbound(2^{j-i}A_{2^j}\mathscr{P}_{2^j}
      \mathscr{Q}_{2^i}^*:i<j)
  =\rbound(A_{2^j}\mathscr{P}_{2^i}
      \mathscr{Q}_{2^j}^*:i<j)\lesssim 1.
\end{split}\end{equation*}

It remains to estimate $A_{2^j}\mathscr{Q}_{2^i}^*$ for $i<j$. We divide this
task into the countable number of cases where $k=j-i\in\Z_+$ is
fixed, aiming to establish sufficiently good $R$-bounds to be able
to sum them up. We start the estimation by writing
\begin{equation}\label{eq:APminusAstart}\begin{split}
  \Big(\Exp &\BNorm{\sum_j\radem_j A_{2^j}
    \mathscr{Q}^{*}_{2^{j-k}}u_j}{L^p(\R^n;X)}^p
    \Big)^{1/p} \\
  &=\Big(\Exp\BNorm{\sum_j\radem_j\sum_{Q\in\triangle_{2^j}}
    1_Q\intav_{Q}\mathscr{Q}^*_{2^{j-k}}u_j}{L^p(\R^n;X)}^p\Big)^{1/p}\\
  &\leq\sum_{\ell=0}^k\Big(\Exp\BNorm{\sum_{j\equiv\ell\,
    \operatorname{mod}\,(k+1)}\radem_j\sum_{Q\in\triangle_{2^j}}
    1_Q\intav_{Q}\mathscr{Q}^*_{2^{j-k}}u_j}{L^p(\R^n;X)}^p\Big)^{1/p}.
\end{split}\end{equation}

We next decompose each of the cubes $Q\in\triangle$ into $2^{k-1}$ parts inductively as follows. Denoting
\begin{equation*}
  \partial_{\delta}E:=\{x\in E:d(x,E^c)\leq\delta\},
\end{equation*}
we set
\begin{equation*}
  Q^1:=\partial_{2^{-k}\ell(Q)}Q,\quad
  Q^m:=\partial_{2^{-k}\ell(Q)}[Q\setminus\bigcup_{\nu=1}^{m-1}Q^{\nu}],
   \quad m=2,\ldots,2^{k-2}.
\end{equation*}
Then $Q^m$ is a union (up to boundaries) of certain dyadic cubes $R\in\triangle_{2^{-k}\ell(Q)}$, and $\abs{Q^m}\geq\abs{Q^{m+1}}$ for all $m<2^{k-1}$. This is  preparation for the application of Lemma~\ref{lem:martingale}  later on.

The right-hand side of \eqref{eq:APminusAstart} may now be rewritten as
\begin{equation*}
  \sum_{\ell=0}^k\Big(\Exp_{\radem'}\Exp_{\radem}
   \BNorm{\sum_{m=1}^{2^{k-1}}\radem_m'
    \sum_{j\equiv\ell}\radem_j
    \sum_{Q\in\triangle_{2^j}}
    1_{Q^m}\intav_Q \mathscr{Q}^*_{2^{j-k}}u_j}{L^p(\R^n;X)}^p\Big)^{1/p},
\end{equation*}
where the randomized $j$ sum, as in Lemma~\ref{lem:signInvar}, does not ``see'' the introduction
of the additional random factors~$\radem_m'$. The UMD space $X$, and then also the Bochner space of functions with values in this space, has some non-trivial Rademacher-type $t>1$, which gives the estimate
\begin{equation*}
  \lesssim\sum_{\ell=0}^k\Big\{\sum_{m=1}^{2^{k-1}}
  \Big(\Exp_{\radem}\BNorm{\sum_{j\equiv\ell}\radem_j
    \sum_{Q\in\triangle_{2^j}}1_{Q^m}\intav_Q \mathscr{Q}^*_{2^{j-k}}u_j}{L^p(\R^n;X)}^p\Big)^{t/p}
  \Big\}^{1/t}.
\end{equation*}

We are now in a position to apply Lemma~\ref{lem:martingale}. For all $m=2, \ldots, 2^{k-1}$, the sets $E(Q)=Q^1$ and $F(Q)=Q^m$ satisfy the assumptions of that Lemma, which means that the summand with $m=1$ above dominates any one of the other summands with $m=2,\ldots,2^{k-1}$. Hence, recalling that $Q^1=\partial_{2^{-k}\ell(Q)}Q$, we may continue with

\begin{equation}\label{used-UMD}
  \lesssim\sum_{\ell=0}^k 2^{k/t}
  \Big(\Exp\BNorm{\sum_{j\equiv\ell}\radem_j
    \sum_{Q\in\triangle_{2^j}}1_{\partial_{2^{j-k}}Q}
    \intav_Q \mathscr{Q}^*_{2^{j-k}}u_j}{L^p(\R^n;X)}^p\Big)^{1/p}.
\end{equation}

Finally, we start making use of the properties of the operators
$\mathscr{Q}^*_t$. For each $Q\in\triangle_{2^j}$, let $\eta_Q\in
C_0^{\infty}(Q)$ be a function with $\eta_Q=1$ in
$Q\setminus\partial_{2^{j-k}}Q$ and $\abs{\nabla\eta_Q}\lesssim
2^{k-j}$. We have
\begin{equation*}\begin{split}
  \int_Q \mathscr{Q}^*_{2^{j-k}}u_j
  &=\int_Q \eta_Q 2^{j-k}(-\diverg)\mathscr{P}_{2^{j-k}}u_j
  +\int_Q (1-\eta_Q)\mathscr{Q}^*_{2^{j-k}}u_j\\
  &= 2^{j-k}\int_Q [\eta_Q,(-\diverg)]\mathscr{P}_{2^{j-k}}u_j
  +\int_Q (1-\eta_Q)\mathscr{Q}^*_{2^{j-k}}u_j,
\end{split}\end{equation*}
where we used the fact that the integral of the divergence of $\eta_Q \mathscr{P}_{2^{j-k}}u_j$ vanishes. We may further observe that $[\eta_Q,(-\diverg)]v=\nabla\eta_Q\cdot v$, and both $\nabla\eta_Q$ and $1_Q-\eta_Q$ are supported on $\partial_{2^{j-k}}Q$, so that both integrals above may be reduced to this smaller set. Thus
\begin{equation}\label{eq:reduceIntegrals}\begin{split}
  &\Big(\Exp\BNorm{\sum_{j\equiv\ell}\radem_j
    \sum_{Q\in\triangle_{2^j}}1_{\partial_{2^{j-k}}Q}
    \intav_Q \mathscr{Q}^*_{2^{j-k}}u_j}{L^p(\R^n;X)}^p\Big)^{1/p}\\
  &=\Big(\Exp\BNorm{\sum_{j\equiv\ell}\radem_j
    \sum_{Q\in\triangle_{2^j}}\frac{\abs{\partial_{2^{j-k}}Q}}{
      \abs{Q}}1_{\partial_{2^{j-k}}Q}\times \\
  &\quad\times\intav_{\doo_{2^{j-k}}Q}\big(2^{j-k}\nabla\eta_Q\cdot
     \mathscr{P}_{2^{j-k}}u_j
  +(1-\eta_Q)\mathscr{Q}^*_{2^{j-k}}u_j\big)}{L^p(\R^n;X)}^p\Big)^{1/p}.
\end{split}\end{equation}
The factors $\abs{\partial_{2^{j-k}}Q}/\abs{Q}$ are equal to $1-(1-2^{1-k})^n\lesssim 2^{-k}$ and may be
extracted outside the summation and the norm. Then we are left
with an expression involving the conditional expectation
projections related to the filtration
\begin{equation*}
  (\sigma(\partial_{2^{j-k}}Q,Q\setminus\partial_{2^{j-k}}Q:
  Q\in\triangle_{2^j}))_{j\equiv\ell\,\operatorname{mod}\, k+1}.
\end{equation*}
These are $R$-bounded under
the UMD assumption, and hence the quantity in \eqref{eq:reduceIntegrals} is majorized by
\begin{equation*}\begin{split}
  &\lesssim 2^{-k}\Big(\Exp\Big\|\sum_{j\equiv\ell}\radem_j
    \sum_{Q\in\triangle_{2^j}}\big(2^{j-k}\nabla\eta_Q\cdot
    \mathscr{P}_{2^{j-k}}u_j+ \\
  & \phantom{ \lesssim 2^{-k}\Big(\Exp\Big\|\sum_{j\equiv\ell}\radem_j\sum_{Q\in\triangle_{2^j}}\big( }  
     +(1_Q-\eta_Q)\mathscr{Q}^*_{2^{j-k}}u_j\big)\Big\|_{L^p(\R^n;X)}^p\Big)^{1/p}\\
  &\lesssim 2^{-k}\Big(\Exp\BNorm{\sum_{j\equiv\ell}\radem_j
    \mathscr{P}_{2^{j-k}}u_j}{L^p(\R^n;X^n)}^p\Big)^{1/p} \\
  &\qquad +2^{-k}\Big(\Exp\BNorm{\sum_{j\equiv\ell}\radem_j
    \mathscr{Q}^*_{2^{j-k}}u_j}{L^p(\R^n;X)}^p\Big)^{1/p},
\end{split}\end{equation*}
where the last estimate used the contraction principle  \ref{contraction-principle}  and
$2^{j-k}\abs{\nabla\eta_Q}\lesssim 1$. Using the $R$-boundedness
of $\mathscr{P}_{2^{j-k}}$ and $\mathscr{Q}^*_{2^{j-k}}$, and substituting back to
\eqref{used-UMD}, we have shown that
\begin{equation*}\begin{split}
  \Big(\Exp &\BNorm{\sum_j\radem_j A_{2^j}\mathscr{Q}^*_{2^{j-k}}u_j}{
   L^p(\R^n;X)}^p\Big)^{1/p} \\
  &\lesssim\sum_{\ell=0}^k 2^{k/t}2^{-k}\Big(\Exp
    \BNorm{\sum_j\radem_j u_j}{L^p(\R^n;X^n)}^p\Big)^{1/p}\\
  &=(k+1)2^{-k/t'}\Big(\Exp\BNorm{\sum_j\radem_j
  u_j}{L^p(\R^n;X^n)}^p\Big)^{1/p}.
\end{split}\end{equation*}

This says that \(\rbound(A_{2^j}\mathscr{Q}^*_{2^{j-k}}:j\in\Z)\lesssim
(k+1)2^{-k/t'}\), and allows us to estimate
\begin{equation*}\begin{split}
  \rbound(2^{\abs{i-j}/2t'}A_{2^j}\mathscr{Q}^*_{2^i}: \; i,j \in \Z \;,\; i<j)
  &\leq\sum_{k=1}^{\infty}\rbound(2^{k/2t'}
    A_{2^j}\mathscr{Q}^*_{2^{j-k}}:j\in\Z)\\
  &\lesssim\sum_{k=1}^{\infty}(k+1) 2^{-k/2t'}
  \lesssim 1.
\end{split}\end{equation*}
We have proved the required $R$-boundedness \eqref{Rbd-for-Schur}
with $\delta=1/2t'=\frac{1}{2}(1-1/t)>0$, where $t>1$ is a
Rademacher-type for $L^p(\R^n;X)$.
\end{proof}

We conclude this section with the following result, which combines most of the estimates achieved so far. Although we will not make direct use of this inequality, but rather the various individual results above, Corollary~\ref{PisA} appears worth recording for the potential further applications of the transference between the dyadic and the harmonic estimates, which it provides.

\begin{corollary}\label{PisA}
Let \(X\) be a UMD space, and \(1<p<\infty\). For \(u \in L^{p}(\R^{n};X)\), we have
\begin{equation*}
  \Exp\BNorm{\sum_k\radem_k(A_{2^k}-\mathscr{P}_{2^k})u}{L^p(\R^n;X)}
  \lesssim \Norm{u}{L^p(\R^n;X)}.
\end{equation*}
For $u\in\range(\Pi)$, the same is true with $P_{2^k}$ in place of
$\mathscr{P}_{2^k}$.
\end{corollary}

\begin{proof}
The first claim is immediate from Propositions~\ref{prop:APminusP}
and \ref{prop:APminusA}, and the second follows
from~\eqref{eq:fcEven}.
\end{proof}


\section{A quadratic $T(1)$ theorem}\label{sect:reduction}

In this section we show that the proof of certain quadratic estimates can be reduced to similar inequalities for the ``principal part'' of the operators involved. 
This will then be applied to our particular operators $Q^B_{2^k}$, and is an analogue of Sections 5.1 and 5.2 in \cite{AKM}. However, we start with the description of a more general situation.

Let $\mathscr{T}=(T_{2^k})_{k\in\Z}$ be an $R$-bounded sequence of
linear operators on $L^p(\R^n;Y)$, where \(1<p<\infty\) and Y is a Banach space, and let $\mathcal{Z}\subseteq L^p(\R^n;Y)$
be a subspace. We say that $\mathscr{T}$ satisfies a
\emph{high-frequency estimate} on $\mathcal{Z}$ if
\begin{equation}\label{eq:highFreqBound}
  \Exp\BNorm{\sum_k\radem_k T_{2^k}(I-\mathscr{P}_{2^k})u}{
    L^p(\R^n;Y)}\lesssim\Norm{u}{L^p(\R^n;Y)}
\end{equation}
for all $u\in\mathcal{Z}$. Concerning the name, note that the symbol of
$I-\mathscr{P}_{2^k}$ is
$(2^{k}\abs{\xi})^2\big(1+(2^{k}\abs{\xi})^2\big)^{-1}$,  which can be thought of as a smooth approximation of the characteristic function of $\{\xi\in\R^n:\abs{\xi}>2^{-k}\}$.

We say that $\mathscr{T}$ satisfies \emph{off-diagonal $R$-bounds}
if the following inequality holds for every $M\in\N$, with the
implied constant only depending on $M$: Whenever $E_k,\
F_k\subset\R^n$ are Borel subsets, $u_k\in L^p(\R^n;Y)$, and
$(t_k)_{k\in\Z}\subseteq\{2^k\}_{k\in\Z}$ are numbers so that
$\dist(E_k,F_k)/t_k>\varrho$ for some $\varrho>0$ and all $k\in\Z$,
there holds
\begin{equation}\label{eq:offDiagonal}\begin{split}
  \Exp &\BNorm{\sum_k\radem_k1_{E_k}
     T_{t_k}1_{F_k}u_k}{L^p(\R^n;Y)} \\
  &\lesssim(1+\varrho)^{-M}\Exp\BNorm{\sum_k\radem_k1_{F_k}u_k}{
     L^p(\R^n;Y)}.
\end{split}\end{equation}
Note that the case $M=0$ follows automatically from the assumed $R$-boundedness of the $T_{2^k}$ and the contraction principle  \ref{contraction-principle} .

Finally, the \emph{principal part} of the operator $T_{2^k}$ is the
operator-valued function $\gamma_{2^k}:\R^n\to\bddlin(Y)$ defined by
(intuitively, ``$\gamma_{2^k}:=T_{2^k}(1)$'')
\begin{equation}\label{eq:prPart}
  \gamma_{2^k}(x)w:=T_{2^k}(w)(x)
  :=\sum_{Q\in\triangle_{2^k}}T_{2^k}(w1_{Q})(x),
  \quad x\in\R^n,\ w\in Y.
\end{equation}
Note that \eqref{eq:offDiagonal} implies that the right-hand side of \eqref{eq:prPart} converges absolutely in $L^p_{\local}(\R^n;Y)$, and this series defines the action of $T_{2^k}$ on the constant function $w$, which lies outside its original domain of definition, namely $L^p(\R^n;Y)$.

We are going to prove the following ``quadratic $T(1)$ theorem'':

\begin{theorem}\label{thm:T1}
Let $Y$ be a UMD space, and $1<p<\infty$. Let the $R$-bounded
operator-sequence $\mathscr{T}=(T_{2^k})_{k\in\Z}$ in
$\bddlin(L^p(\R^n;Y))$ satisfy the high-frequency
estimate~\eqref{eq:highFreqBound} on a subspace $\mathcal{Z}\subseteq
L^p(\R^n;Y)$, and the off-diagonal $R$-bounds~\eqref{eq:offDiagonal}.
Then there holds
\begin{equation*}
  \Exp\BNorm{\sum_k\radem_k \big(T_{2^k}-\gamma_{2^k}A_{2^k}\big)u}{
     L^p(\R^n;Y)}
  \lesssim\Norm{u}{L^p(\R^n;Y)},\qquad
   u\in \mathcal{Z}.
\end{equation*}
Thus $\mathscr{T}$ satisfies a quadratic estimate on $\mathcal{Z}$ if and only
if its principal part does.
\end{theorem}

Before going into the proof, let us indicate the consequences for our primary case of interest,
 which is the vector-valued analogue of 
Proposition 5.5 in \cite{AKM}:

\begin{theorem}\label{reduction}
Let $X$ be a UMD Banach space, $1<p<\infty$, and $\Pi_{B}$ be an
$R$-bisectorial perturbed Hodge-Dirac operator on $L^p(\R^n;X^N)$.
Let $\gamma_{2^k}$ denote the principal part of $Q_{2^k}^B$. Then
there holds:
\begin{equation*}
  \Exp\BNorm{\sum_{k\in\Z}\radem_k(Q_{2^k}^B-\gamma_{2^k}A_{2^k})u}{L^p(\R^n;X^N)}
  \lesssim\Norm{u}{L^p(\R^n,X^N)},\quad
  \forall u\in\range(\Gamma),
\end{equation*}
and the operators $\gamma_{2^k}(x)$ are multiplications
by complex $N\times N$-matrices.
\end{theorem}

The quadratic estimate is obviously implied by Theorem~\ref{thm:T1} as
soon as we check that $(Q^B_{2^k})_{k\in\Z}$ satisfies the
high-frequency estimate on $\range(\Gamma)$ and the off-diagonal
$R$-bounds. This is the content of the next two results below. The form of the principal part follows readily from the
definition~\eqref{eq:prPart} and the fact that the operators $Q^B_{2^k}$ on $L^p(\R^n;X^N)$ are tensor extensions of operators on
$L^p(\R^n;\C^N)$.

\begin{lemma}\label{lem:QminusQP}
The family $(Q^B_{2^k})_{k\in\Z}$ satisfies the high-frequency
estimate \eqref{eq:highFreqBound} on $\range(\Gamma)\subset
L^p(\R^n;X^N)$.
\end{lemma}

\begin{proof}
It follows from~\eqref{eq:fcEven} that $\mathscr{P}_{2^k}u= P_{2^k}u$
for $u\in\range(\Gamma)$, so it suffices to prove the modified claim
with $P_{2^k}$ in place of $\mathscr{P}_{2^k}$.

Let \(\mathbb{P}^{1}\) denote the projection of
\begin{equation*}
  L^p(\R^n;X^N)
  =L^p(\R^n;X^{n_1})\oplus L^p(\R^n;X^{n_2})
\end{equation*}
onto \(L^p(\R^n;X^{n_1})\). Since \(u \in\range(\Gamma)\), a
straightforward manipulation using the structure of the operators
shows that
\begin{equation*}
  Q_{t}^{B}(I-P_{t})u = Q_{t}^{B}t\Gamma Q_{t}u =
  (I-P_{t}^{B})\mathbb{P}^{1}Q_{t}u.
\end{equation*}
Since \(\{(I-P_{t}^{B})\mathbb{P}^{1} \;;\; t \geq 0\}\) is
\(R\)-bounded, this gives
\begin{equation*}\begin{split}
  \Exp\BNorm{\sum_k\radem_k Q^B_{2^k}(I-P_{2^k})u}{L^p(\R^n;X^N)}
  &\lesssim\Exp\BNorm{\sum_k\radem_k Q_{2^k}u}{L^p(\R^n;X^N)} \\
  &\lesssim\Norm{u}{L^p(\R^n;X^N)},
\end{split}\end{equation*}
where the last inequality follows from Proposition \ref{funcalc}.
\end{proof}

 The following Proposition is the vector-valued analogue of 
Proposition 5.2 in \cite{AKM}. 

\begin{proposition}\label{offdiag}
The family $(Q_{2^k}^B)_{k\in\Z}$ satisfies the off-diagonal
$R$-bounds \eqref{eq:offDiagonal}.
\end{proposition}

\begin{proof}
It is sufficient to prove this result for \(R_{t_{k}}^{B}\) instead
of \(Q_{t_{k}}^{B}\) since \(Q_{t_{k}}^{B} =
\frac{i}{2}(R_{t_{k}}^{B} - R_{-t_{k}}^{B})\). We proceed by
induction on \(M\). The case \(M=0\) follows from Kahane's
contraction principle  \ref{contraction-principle}  and the R-bisectoriality of \(\Pi_{B}\).  Now
assume it is true for some \(M\geq 0\), and consider
\begin{equation*}
  \tilde{E}_{k} = \{x\in\R^{n} \;;\; \dist(x,E_{k}) < \frac{1}{2}\dist(x,F_{k})\}
\end{equation*}
and \(\eta_{k}\) a cutoff function supported in \(\tilde{E}_{k}\)
with \((\eta_{k})_{|E_{k}}=1\) and \(\|\nabla\eta_{k}\|_{\infty}
\leq 4/\dist(E_{k},F_{k})\). Denoting by \([T,S]=TS-ST\) the
commutator of two operators we have
\begin{equation*}
  [\eta_{k}I,R_{t_{k}}^{B}] = it_{k}R_{t_{k}}^{B}([\Gamma,\eta_{k}I]+B_{1}[\Gamma^{*},\eta_{k}I]B_{2})R_{t_{k}}^{B}.
\end{equation*}
Using $R$-bisectoriality, and the fact that
\([\Gamma,\eta_{k}I]+B_{1}[\Gamma^{*},\eta_{k}I]B_{2}\) is a
multiplication by an \(L^{\infty}\) function bounded by \(\|\nabla
\eta_{k}\|_{\infty}\), we thus have
\begin{equation*}\begin{split}
  \Exp &\BNorm{\sum_k \varepsilon_{k}1_{E_{k}}R^{B} _{t_{k}}1_{F_{k}}u_{k}}{} \\
  &\lesssim\Exp\BNorm{\sum_k \varepsilon_{k}[\eta_{k}I,R^{B} _{t_{k}}]1_{F_{k}}u_{k}}{} \\
  &\lesssim \Exp\BNorm{\sum_k \varepsilon_{k}  it_{k}R_{t_{k}}^{B}
    \big([\Gamma,\eta_{k}I]+B_{1}[\Gamma^{*},\eta_{k}I]B_{2}\big)
    1_{\tilde{E}_{k}}R_{t_{k}}^{B}1_{F_{k}}u_{k}}{} \\
  &\lesssim \sup_{j \in \Z}\abs{t_{j}}\Norm{\nabla\eta_{j}}{\infty}
    \Exp\BNorm{\sum_k \varepsilon_{k}
     1_{\tilde{E}_{k}}R_{t_{k}}^{B}1_{F_{k}}u_{k}}{} \\
  &\lesssim \frac{1}{\rho}\,\Exp\BNorm{\sum_k \varepsilon_{k}
   1_{\tilde{E}_{k}}R_{t_{k}}^{B}1_{F_{k}}u_{k}}{},
\end{split}\end{equation*}
and we may apply the induction assumption to the remaining quantity.
\end{proof}

This completes the proof that Theorem \ref{reduction} is a consequence of Theorem \ref{thm:T1}. 
We now return to the Quadratic $T(1)$ Theorem~\ref{thm:T1}. In
proving this result, we decompose
\begin{equation*}
  T_t-\gamma_t A_t
  =T_t(I-\mathscr{P}_t)+(T_t-\gamma_t A_t)\mathscr{P}_t+\gamma_t A_t(\mathscr{P}_t-I),
\end{equation*}
where the different summands on the right will be analyzed
separately. The first one, of course, is immediately handled by the
assumed high-frequency estimate.

\begin{lemma}\label{lem:gammaRbdd}
Under the assumptions of Theorem~\ref{thm:T1}, the principal part
operators \((\gamma_{2^{k}}A_{2^{k}})_{k \in \Z}\) are $R$-bounded
on $L^p(\R^n;X^N)$.
\end{lemma}

\begin{proof}
For \((u_{k})_{k\in\Z} \subset L^p(\R^n;X^N)\) we have
\begin{equation*}\begin{split}
  \Exp &\BNorm{\sum_{k\in\Z}\radem_k\gamma_{2^k}A_{2^k}u_k}{}
  =\Exp\BNorm{\sum_{k\in\Z}\radem_k\sum_{Q\in\triangle_{2^k}}1_Q
    T_{2^k}\ave{u_k}_Q}{} \\
  &\leq\sum_{m\in\Z^n}\Exp\BNorm{\sum_{k\in\Z}\radem_k\sum_{Q\in\triangle_{2^k}}
    1_Q T_{2^k}(1_{Q+2^k m}\ave{u_k}_Q)}{} \\
  &\lesssim\sum_{m\in\Z^n}(1+\abs{m})^{-M}
   \Exp\BNorm{\sum_{k\in\Z}\radem_k\sum_{Q\in\triangle_{2^k}}
    1_{Q+2^k m}\ave{u_k}_Q}{} \\
  &\lesssim\sum_{m\in\Z^n}(1+\abs{m})^{-M}
   \log(2+\abs{m})\,\Exp\BNorm{\sum_{k\in\Z}\radem_k u_k}{},
\end{split}\end{equation*}
where the last two estimates where applications of the off-diagonal
estimates (and sign-invariance), and Corollary~\ref{cor:Figiel},
respectively. The series is summable for $M>n$.
\end{proof}

 The next Lemma is the vector-valued analogue of 
Proposition~5.5 in~\cite{AKM}. 

\begin{lemma}\label{lem:multiplierPart}
Under the assumptions of Theorem~\ref{thm:T1}, for all $u\in
L^p(\R^n;Y)$ there holds
\begin{equation*}
  \Exp\BNorm{\sum_{k\in\Z}\radem_k(T_{2^k}-\gamma_{2^k}A_{2^k})
    \mathscr{P}_{2^k}u}{L^p(\R^n,Y)}
  \lesssim\Norm{u}{L^p(\R^n,Y)}.
\end{equation*}
\end{lemma}

\begin{proof}
We first observe that it suffices to prove a modified assertion with
$\mathscr{P}_{2^k}$ replaced by $\mathscr{P}_{2^k}^M$. Indeed, this
follows at once from the $R$-boundedness of $T_{2^k}$ and
$\gamma_{2^k}A_{2^k}$ combined with Lemma~\ref{PtVsPtN}.

As for the new claim, denote $v_k:=\mathscr{P}_{2^k}^M u$. Then
\begin{equation}\label{eq:multiplierPart}\begin{split}
  \Exp &\BNorm{\sum_k\radem_k(T_{2^k}-\gamma_{2^k}A_{2^k})v_k}{} \\
  &=\Exp\BNorm{\sum_k\radem_k\sum_{Q\in\triangle_{2^k}}1_Q
     T_{2^k}\big(v_k-\ave{v_k}_Q)}{} \\
  &\leq\sum_{m\in\Z^n}\Exp\BNorm{\sum_k\radem_k\sum_{Q\in\triangle_{2^k}}
     1_Q T_{2^k}\big(1_{Q-2^k m}(v_k-\ave{v_k}_Q)\big)}{} \\
  &\lesssim\sum_{m\in\Z^n}(1+\abs{m})^{-M}\,
   \Exp\BNorm{\sum_k\radem_k\sum_{Q\in\triangle_{2^k}}
     1_Q(v_k-\ave{v_k}_{Q+2^k m})}{}
\end{split}\end{equation}
where we used the off-diagonal estimates. By the Poincar\'e inequality
(Proposition~\ref{poincare}) and Lemma \ref{lem:multiplier}, the last factor is majorized by
\begin{equation*}\begin{split}
  \int_{[-1,1]^n}\int_0^1
    &\Exp\BNorm{\sum_k\radem_k 2^k(m+z)\cdot\nabla\,\tau_{t2^k(m+z)}
     \mathscr{P}_{2^k}^M u}{}
    \ud t \ud z \\
  &\lesssim(1+\abs{m})^{n+1} \Norm{u}{}.
\end{split}\end{equation*}
Substituing this back to \eqref{eq:multiplierPart}, we find that the series sums up to \(\lesssim\Norm{u}{}\) provided that we choose \(M>2n+1\).
\end{proof}

\begin{proof}[of Theorem~\ref{thm:T1}]
We have
\begin{equation*}\begin{split}
  \Exp &\BNorm{\sum_k\radem_k(T_{2^k}-\gamma_{2^k}A_{2^k})u}{}
  \lesssim \Exp\BNorm{\sum_k\radem_k T_{2^k}(I-\mathscr{P}_{2^k})u}{} \\
   &+\Exp\BNorm{\sum_k\radem_k(T_{2^k}-\gamma_{2^k}A_{2^k})\mathscr{P}_{2^k}u}{}
   +\Exp\BNorm{\sum_k\radem_k\gamma_{2^k}A_{2^k}(\mathscr{P}_{2^k}-I)u}{}.
\end{split}\end{equation*}
For $u\in Z\subset L^p(\R^n;Y)$, the upper bound $\Norm{u}{}$ for
the first term follows from the assumed high-frequency estimate, for
the second term from Lemma~\ref{lem:multiplierPart}, and for the
third one from Lemma~\ref{lem:gammaRbdd} and
Proposition~\ref{prop:APminusA} together with the observation that
$A_{2^k}=A_{2^k}A_{2^k}$.
\end{proof}

In order to estimate the principal term
\begin{equation}\label{eq:principal}
  \Exp\BNorm{\sum_{k\in\Z}\radem_k
    \gamma_{2^{k}}A_{2^{k}}u}{L^p(\R^n;X^N)},
  \qquad u\in\range(\Gamma),
\end{equation}
we need a version of Carleson's inequality. This is achieved in
Section~\ref{carleson} by using the Rademacher maximal function,
which we next study.


\section{The Rademacher maximal function}\label{maximal}

We recall the definition of the Rademacher maximal function, here stated in an equivalent but slightly different way from Section~\ref{prelim}:
\begin{equation*}\begin{split}
  M_R u(x):=\sup\Big\{
    \Exp &\BNorm{\sum_{Q\owns x}\radem_Q\lambda_Q\ave{u}_Q}{X}:\\
    &(\lambda_Q)_{Q\in\triangle}\text{
    finitely non-zero with }\sum_{Q\in\triangle}\abs{\lambda_Q}^2\leq 1\Big\}.
\end{split}\end{equation*}
We will also find it convenient to consider the following linearized version:
\begin{equation*}
  \mathscr{M}_Ru(x):\ell^2(\triangle)\to\Rad(X),
  (\lambda_Q)_{Q\in\triangle}
  \mapsto\sum_{Q\owns x}\radem_Q\lambda_Q\ave{u}_Q,
\end{equation*}
which satisfies \(M_Ru(x)=\Norm{\mathscr{M}_R u(x)}{\bddlin(\ell^2,\Rad(X))}\).

The RMF property of a Banach space \(X\) was defined in terms of the \(L^2\)-boundedness of \(M_R\), but the next result shows that the exponent \(2\) is not relevant:

\begin{proposition}\label{prop:allP}
Let $X$ be a Banach space, and consider the assertion
\begin{equation}\label{eq:MRbdd}
  M_R:L^p(\R^n;X)\to L^p(\R^n)
  \text{ is bounded}.
\end{equation}
If \eqref{eq:MRbdd} is true for one $p\in(1,\infty)$,
then it is true for all $p\in(1,\infty)$.
\end{proposition}

\begin{proof}
 It suffices to prove the same for the equivalent statement
\begin{equation}\label{eq:MRbddlin}
  \mathscr{M}_R:L^{p}(\R^n,X)\to L^{p}(\R^n,\bddlin(\ell^2,\Rad(X)))
  \text{ is bounded}.
\end{equation}
Suppose that \eqref{eq:MRbddlin} is true for some
$p\in(1,\infty)$. Let $a$ be a dyadic atom of
$H^1(\R^n,X)$, i.e., $\supp a\subseteq Q$, a dyadic cube,
$\Norm{a}{\infty}\leq\abs{Q}^{-1}$ and $\int a(x)\ud x=0$. Then
$\ave{a}_{Q'}\neq 0$ only if $Q'\subset Q$. Hence
\begin{equation*}\begin{split}
  &\Norm{\mathscr{M}_R u}{L^1(\R^n,\bddlin(\ell^2,\Rad(X)))}
  =\Norm{\mathscr{M}_R u}{L^1(Q,\bddlin(\ell^2,\Rad(X)))} \\
  &\leq\abs{Q}^{1/p'}\Norm{\mathscr{M}_R u}{L^p(\R^n,\bddlin(\ell^2,\Rad(X)))} \\
  &\lesssim\abs{Q}^{1/p'}\Norm{u}{L^p(\R^n,X)}
  \leq\abs{Q}^{1/p'}\abs{Q}^{1/p}\Norm{u}{\infty}\leq 1.
\end{split}\end{equation*}
It follows that $\mathscr{M}_R:H^1(\R^n,X)\to
L^1(\R^n,\bddlin(\ell^2,\Rad(X)))$ boundedly.

Let then $u\in L^{\infty}(\R^n,X)$ and let $Q$ be a dyadic
cube. It is easy to see that
\begin{equation*}
  1_Q[\mathscr{M}_Ru-\ave{\mathscr{M}_Ru}_Q]
  =\mathscr{M}_R(1_Q[u-\ave{u}_Q]).
\end{equation*}
It follows that
\begin{equation*}\begin{split}
  &\Norm{\mathscr{M}_R u}{\BMO(\R^n,\bddlin(\ell^2,\Rad(X)))} \\
  &=\sup_{Q\in\triangle}\frac{1}{\abs{Q}}
    \Norm{\mathscr{M}_Ru-\ave{\mathscr{M}_Ru}_Q}{
    L^1(Q,\bddlin(\ell^2,\Rad(X)))}\\
  &=\sup_{Q\in\triangle}
    \Norm{\mathscr{M}_R(\abs{Q}^{-1}1_Q[u-\ave{u}_Q]}{
    L^1(\R^n,\bddlin(\ell^2,\Rad(X)))}.\\
\end{split}\end{equation*}
But $\abs{Q}^{-1}1_Q[u-\ave{u}_Q]$ is $2\Norm{u}{\infty}$ times an
atom of $H^1(\R^n,X)$. Hence, by what we already showed, we also
find that
$\mathscr{M}_R:L^{\infty}(\R^n,X)\to\BMO(\R^n,\bddlin(\ell^2,\Rad(X)))$
boundedly. Now interpolation gives the assertion.
\end{proof}

\begin{remark}
Given a dyadic cube \(Q\in\triangle\), it also makes sense to consider \(M_R\) as an operator acting in \(L^p(Q;X)\). In this case one may restrict the summation in the definition to
\begin{equation*}
  \sum_{R:x\in R\subseteq Q}\radem_R\lambda_R\ave{u}_R.
\end{equation*}
An obvious restriction argument now shows that \(M_R:L^p(Q;X)\to L^p(Q)\), with the norm independent of \(Q\), if \(X\) has RMF.
\end{remark}

We do not yet fully understand how the RMF property relates to established Banach space notions. Since we need to assume this kind of  inequality to be able to carry out the estimates in the subsequent sections, we next provide some sufficient conditions, which imply this property. In Appendix~\ref{app:l1} we also give a counterexample to show that RMF is indeed a nontrivial property not shared by every Banach space; more precisely, it fails in the sequence space \(\ell^1\). Our first sufficient condition, Rademacher type \(2\), is the easiest one, but not very useful for our applications, since this condition is not self-dual and the condition that both \(X\) and \(X^*\) have type \(2\) is very restrictive, indeed, equivalent to \(X\) being isomorphic to a Hilbert space. On the other hand, the other two classes of spaces with RMF --- UMD function lattices and reflexive noncommutative \(L^p\) spaces --- are both self-dual, and they cover the most important concrete examples of UMD spaces.

\subsubsection*{Spaces of type \(2\).}

If $X$ has type $2$, then $M_R u(x)\lesssim Mu(x)$,
where $M$ is the usual dyadic maximal function. In fact,
\begin{equation}\label{eq:MRvsM}
  \Exp\BNorm{\sum_{k}\radem_k\lambda_kA_{2^k}u(x)}{X}
  \lesssim\Big(\sum_{k}\abs{\lambda_k}^2
    \Norm{A_{2^k}u(x)}{X}^2\Big)^{1/2}
\end{equation}
in this case, and the supremum over
$\Norm{\lambda}{\ell^2(\Z)}\leq 1$ of the right-hand side
is $\sup_{k}\norm{A_{2^k}u(x)}{X}=Mu(x)$.

\begin{remark} If $X$ has cotype
$2$, then the reverse estimate holds in \eqref{eq:MRvsM}, and hence
$M_R u(x)\gtrsim Mu(x)$. Thus $M_R u(x)\eqsim Mu(x)$ if $X$ is
(isomorphic to) a Hilbert space.
\end{remark}

\begin{remark}
In \cite{james}, James constructed a non reflexive Banach space with type 2 (and thus with the RMF property).
This means, in particular, that RMF does not imply UMD.
\end{remark}

\subsubsection*{UMD function lattices.}

Suppose now that $X$ is a Banach lattice of (equivalence classes of) measurable functions on some \(\sigma\)-finite measure space \((S,\Sigma,\mu)\). This means that \(X\) is a Banach space of such functions and, in addition,
\begin{itemize}
  \item it contains the pointwise real and imaginary parts of any two functions \(\xi,\eta\in X\), and the pointwise maximum and minimum of any two real function \(\xi,\eta\in X\);
  \item if the pointwise absolute values satisfy \(\abs{\xi}\leq\abs{\eta}\), then \(\Norm{\xi}{X}\leq\Norm{\eta}{X}\).
\end{itemize}
Obvious examples are the \(L^{p}(\mu)\) and spaces of continuous functions; also any Banach space with an unconditional basis may be viewed as a Banach lattice of functions defined on \(\Z_+\). One can also give an abstract definition of a Banach lattice without a postulated function space structure (see e.g. \cite{aliprantis}), but we restrict ourselves to the concrete situation, which is the context where Banach lattices with the UMD property have been studied by Rubio de Francia \cite{Rubio}. In this situation, the harmonic analysis in \(L^{p}(\R^{n};X)\) is much closer to the scalar valued case than on a general UMD space, since one can use square functions similar to their \(L^{p}(\R^{n};\C)\) counterparts, and there is also the following natural notion of a maximal function.

The (dyadic) lattice maximal function $M_{\text{lattice}}$ is defined by
\begin{equation*}
  M_{\text{lattice}}u(x):=\sup_{Q\owns x}\abs{\ave{u}_Q},
\end{equation*}
which is again an $X$-valued function.  Suppose $X$ is UMD (and thus has finite cotype), then
\begin{equation*}\begin{split}
  \Exp\BNorm{\sum_{Q\owns x}\radem_Q\lambda_Q\ave{u}_Q}{X}
  &\lesssim\BNorm{\Big(\sum_{Q\owns x} 
    \abs{\lambda_Q}^2\abs{\ave{u}_Q}^2\Big)^{1/2}}{X} \\
  &\leq\BNorm{\Big(\sum_{Q\owns x}
    \abs{\lambda_Q}^2\Big)^{1/2}
    \sup_{Q\owns x}\abs{\ave{u}_Q}}{X},
\end{split}\end{equation*}
so that we have the domination $M_R u(x)\lesssim
\Norm{M_{\text{lattice}}u(x)}{X}$. By a result of Rubio de Francia \cite{Rubio}, we know that
$\Norm{M_{\text{lattice}}u}{L^p(\mu,X)}\lesssim\Norm{u}{L^p(\mu,X)}$,
and hence $\Norm{M_R u}{L^p(\mu)}\lesssim\Norm{u}{L^p(\mu,X)}$ for
all $1<p<\infty$.\\

\subsubsection*{Noncommutative $L^{p}$ spaces.}
We now turn to the case where \(X\) is a noncommutative \(L^{p}\) space \(L^p(N,\tau)\) on a von Neumann algebra \(N\) with a normal semifinite faithful trace~\(\tau\). In this setting, analogues of many important results from Banach space theory and harmonic analysis have recently been found. See \cite{pisierxu} for the definition, more information and references. We here presuppose a modest knowledge of these notions, and only mention that the \(L^p(N,\tau)\) are spaces of (bounded linear) operators (acting on some Hilbert space),  which generalize the ``commutative'' \(L^{p}(\mu)\) spaces, the trace playing the r\^ole of an integral. The simplest examples, besides \(L^{p}\), are the Schatten ideals \(S^{p}\) of bounded linear operators \(A\) such that \(\operatorname{tr}((A^{*}A)^{p/2})\) is finite, where \(\operatorname{tr}\) denotes the usual trace. The reader who is not interested in the applications of our results in the noncommutative context, may very well jump to the beginning of the next section.

The following ``noncommutative Doob's maximal inequality'' was established by M.~Junge \cite{Junge}:

\begin{theorem}[Junge]\label{thm:Junge}
Let $1<p\leq\infty$ and $u\in L^p(N,\tau)$. Let $(N_i)$ be an
increasing sequence of von Neumann subalgebras of $N$, with
associated conditional expectations $E_i$. Then there exist
$a,b\in L^{2p}(N,\tau)$ and contractions $y_i\in N$ such that
\[
  E_iu=ay_ib,\qquad \Norm{a}{2p}\Norm{b}{2p}\lesssim_p
    \Norm{u}{p}.
\]
\end{theorem}

In particular (cf.\ \cite{Junge}, Remark~5.5),
Theorem~\ref{thm:Junge} applies in the case when
\begin{equation*}
  N=L^{\infty}(\mathscr{F})\bar\otimes M,
\end{equation*}
where $L^{\infty}(\mathscr{F})$ is a usual commutative $L^{\infty}$
space, and $N_i=L^{\infty}(\mathscr{F}_i)\bar\otimes M$ for some
sub-$\sigma$-algberas $\mathscr{F}_i\subset\mathscr{F}$. Then
$L^p(N)\eqsim L^p(\mathscr{F},L^p(M))$ is the Bochner space of
$L^p$ functions with values in the noncommutative space $L^p(M)$,
and $E_i$ are the (tensor extensions of) usual conditional
expectation operators. 
In our case \(E_{i} = A_{2^{i}}\), but the argument is valid for general sequences of conditional expectations.

\begin{corollary}
Let $1<p,q<\infty$, let $X=L^q(M)$ and $u\in L^p(\mathscr{F},X)$.
Then
\begin{equation*}
  \Norm{M_R u}{L^p(\mathscr{F})}
  \lesssim_{p,q}\Norm{u}{L^p(\mathscr{F},X)}.
\end{equation*}
\end{corollary}

\begin{proof}
By Proposition~\ref{prop:allP}, it suffices to prove the case
$p=q$. Then $L^p(\mathscr{F},L^p(M))=L^p(N)$, with \(N=L^{\infty}(\mathscr{F})\bar\otimes M\), is itself a
noncommutative $L^p$ space. By Theorem~\ref{thm:Junge}, there
exist $a,b\in L^{2p}(N)=L^{2p}(\mathscr{F},L^{2p}(M))$ and
contractions $y_j\in N$ such that
\begin{equation}\label{eq:Junge}
  E_j u(x)=a(x)y_j(x)b(x),\qquad
  \Norm{a}{L^{2p}(N)}\Norm{b}{L^{2p}(N)}\lesssim_p
  \Norm{u}{L^p(N)}.
\end{equation}
Then we have, by the noncommutative H\"older inequality,
\begin{equation*}\begin{split}
  \Exp\BNorm{\sum_j\radem_j\lambda_j E_j u(x)}{L^p(M)}
  &=\Exp\BNorm{a(x)\sum_j\radem_j\lambda_j y_j(x)b(x)}{L^p(M)}\\
  &\leq\Exp\Norm{a(x)}{L^{2p}(M)}
    \BNorm{\sum_j\radem_j\lambda_j y_j(x)b(x)}{L^{2p}(M)}.
\end{split}\end{equation*}
Now $2p>2$, so that the space $L^{2p}(M)$ has type $2$. Hence
\begin{equation*}\begin{split}
  \Exp &\BNorm{\sum_j\radem_j\lambda_j y_j(x)b(x)}{L^{2p}(M)}
  \lesssim_p\Big(\sum_j\Norm{\lambda_j y_j(x)b(x)}{
    L^{2p}(M)}^2\Big)^{1/2}\\
  &\leq\Big(\sum_j\abs{\lambda_j}^2\Big)^{1/2}
    \Norm{b(x)}{L^{2p}(M)}
  \leq\Norm{b(x)}{L^{2p}(M)}.
\end{split}\end{equation*}
Combining the previous estimates, we have shown that
\begin{equation*}
  M_Ru(x)\lesssim_p \Norm{a(x)}{L^{2p}(M)}
    \Norm{b(x)}{L^{2p}(M)},
\end{equation*}
and hence, by H\"older's inequality and \eqref{eq:Junge},
\begin{equation*}
  \Norm{M_R u}{L^p(\mathscr{F})}
  \lesssim_p\Norm{a}{L^{2p}(\mathscr{F};L^{2p}(M))}
   \Norm{b}{L^{2p}(\mathscr{F};L^{2p}(M))}
  \lesssim_p\Norm{u}{L^p(\mathscr{F};L^p(M))},
\end{equation*}
which completes the proof.
\end{proof}
The results of this section constitute a proof of Proposition \ref{radmax}.


\section{An $L^{p}$ version of Carleson's inequality}\label{carleson}

We next establish a vector-valued $L^p$ version of Carleson's inequality for Carleson measures. For $p\neq 2$, it appears to be new even in the scalar-valued case. We wish to mention that the proof of this inequality is significantly inspired by the work of N.~H.\ Katz and M.~C.\ Pereyra \cite{KP,Pereyra}, although none of their specific results is explicitly needed.

Let  $b=(b_R)_{R\in\triangle}$ be a finitely non-zero sequence of
measurable scalar-valued functions, such that $\supp b_Q\subseteq Q$. For each \(Q \in \triangle\) we denote 
\begin{equation*}\begin{split}
  \Norm{b}{\Car^p(Q)}&:=
  \sup_{S \in \triangle \,,\, S\subseteq Q}\Big(\frac{1}{\abs{S}}\int_S
    \Exp\Babs{\sum_{R\subset S}\radem_R b_R(x)}^p
    \ud x\Big)^{1/p}\\
  &\eqsim
  \sup_{S \in \triangle \,,\,  S\subseteq Q}\Big(\frac{1}{\abs{S}}\int_S
    \Big[\sum_{R\subset S}\abs{b_R(x)}^2\Big]^{p/2}
    \ud x\Big)^{1/p}.
\end{split}\end{equation*}
Let us write
$\Norm{b}{\Car^p(\R^n)}:=\sup_{Q\in\triangle}\Norm{b}{\Car^p(Q)}$.
For $p=2$, this is just (the square-root of) the Carleson constant of the measure
\begin{equation*}
  \ud\mu(x,t)=\sum_{Q\in\triangle}b_Q(x)1_{\left]\ell(Q)/2,\ell(Q)\right]}(t)
    \ud x\frac{\ud t}{t}.
\end{equation*}

For the moment, fix a cube $Q\in\triangle$, and denote by $\mu$
the normalized Lebesgue measure, $\mu(E):=\abs{E}/\abs{Q}$, on
measurable subsets of $Q$. We recall the definition of Lorentz
spaces $L^{p,q}(\mu,X)$. A measurable function $u:Q\to X$ belongs
to $L^{p,q}(\mu,X)$ if
\begin{equation*}
  \Norm{u}{L^{p,q}(\mu,X)}:=
  \Big(\int_0^{\infty}\big[t\mu(\Norm{u(\cdot)}{X}>t)^{1/p}\big]^q\frac{\ud t}{t}
   \Big)^{1/q}
\end{equation*}
is finite. We are now ready to state:

\begin{lemma}\label{lem:coreOfCarleson}
Let $X$ be a Banach space with type $t\geq 1$, and let
$1\leq p<\infty$. Then
\begin{equation*}\begin{split}
  &\Big(\frac{1}{\abs{Q}}\int_Q\Exp\BNorm{
   \sum_{R \in \triangle \,,\, R\subseteq Q}\radem_R b_R(x)\ave{u}_R}{X}^p
  \Big)^{1/p} \\
  &\lesssim\Norm{b}{\Car^p(Q)}\times\begin{cases}
    \Norm{M_R u}{L^{p}(\mu)} & \text{if }1\leq p\leq t,\\
    \Norm{M_R u}{L^{p,t}(\mu)} & \text{if }t<p<\infty.
  \end{cases}
\end{split}\end{equation*}
\end{lemma}

\begin{proof}
Let us fix some $A>0$ and denote
\begin{equation*}
  \mathscr{G}_k:=\Big\{S\subseteq Q:
  \sup_{\Norm{\lambda}{\ell^2}\leq 1}\Exp\BNorm{\sum_{R:S\subseteq R\subseteq Q}
    \radem_R\lambda_R\ave{u}_R}{X}\leq A\cdot 2^k\Big\}.
\end{equation*}
Let us also denote by $\mathscr{F}_k$ the set of maximal dyadic
cubes $S\subseteq Q$ such that $S\notin\mathscr{G}_k$.

Then every $R\notin\mathscr{G}_k$ satisfies $R\subseteq S$ for a
unique $S\in\mathscr{F}_k$. Moreover,
$\mathscr{G}_k\subseteq\mathscr{G}_{k+1}$, and every $S\subseteq
Q$ belongs to $\mathscr{G}_k$ for a sufficiently large $k$. We
write $\mathscr{Q}_0:=\mathscr{G}_0$ and
$\mathscr{Q}_k:=\mathscr{G}_k\setminus\mathscr{G}_{k-1}$ for
$k=1,2,\ldots$ Then
\begin{equation*}
  \sum_{R\subseteq Q}\radem_R b_R(x)\ave{u}_R
  =\sum_{k=0}^{\infty}\sum_{R\in\mathscr{Q}_k}
   \radem_R b_R(x)\ave{u}_R,
\end{equation*}
and, by sign-invariance, 
\begin{equation*}
\Exp \|\sum_{k=0}^{\infty}\sum_{R\in\mathscr{Q}_k}
   \radem_R b_R(x)\ave{u}_R\| \eqsim
\Exp\Exp'  \|\sum_{k=0}^{\infty}\radem_k'\sum_{R\in\mathscr{Q}_k}
   \radem_R b_R(x)\ave{u}_R\|,
\end{equation*}
where $\radem_k'$ are an independent sequence of Rademacher
variables. Let us denote $q:=\min\{p,t\}$, so that $X$ has type $q$.

Then, by the definition of type,
\begin{equation*}
   \Exp\Exp'\BNorm{\sum_{k=0}^{\infty}\radem_k'\sum_{R\in\mathscr{Q}_k}
   \radem_R b_R(x)\ave{u}_R}{X}^p
  \lesssim\Big(\sum_{k=0}^{\infty}\Exp\BNorm{\sum_{R\in\mathscr{Q}_k}
   \radem_R b_R(x)\ave{u}_R}{X}^q\Big)^{p/q}.
\end{equation*}

Now consider a fixed $x\in Q$. Suppose first that there is a smallest dyadic cube $S$ such that $x\in S\in\mathscr{Q}_k$. Then
\begin{equation}\label{eq:stoppingArg}\begin{split}
  \Exp&\BNorm{\sum_{R\in\mathscr{Q}_k}
   \radem_R b_R(x)\ave{u}_R}{X}^q
  =\Exp\BNorm{\sum_{S\subseteq R\subseteq Q}\radem_R b_R(x)1_{\mathscr{Q}_k}(R)\ave{u}_R}{X}^q \\
  &\lesssim (A2^k)^q\Exp\Babs{\sum_{S\subseteq R\subseteq Q}\radem_R b_R(x)1_{\mathscr{Q}_k}(R)}^q
  =(A2^k)^q\Exp\Babs{\sum_{R\in\mathscr{Q}_k}\radem_R b_R(x)}^q,
\end{split}\end{equation}
where the estimate employed the fact that $S\in\mathscr{Q}_k\subseteq\mathscr{G}_k$, the defining property of $\mathscr{G}_k$ with $\lambda_R=b_R(x)1_{\mathscr{Q}_k}(R)$, and the equivalence of the $\ell^2$ norm and the randomized norm for scalar sequences.

If there is no smallest $S$, then \eqref{eq:stoppingArg} remains true with ``$\lim_{S\downarrow\{x\}}$'' in front of the two intermediate expressions, where $S$ runs through the decreasing sequence of dyadic cubes containing $x$. In either case, the final estimate between the left-hand and the right-hand side is the same.

Substituting back and using 
the triangle inequality in $L^{p/q}(\mu)$, we have
\begin{equation*}\begin{split}
  \Big(&\frac{1}{\abs{Q}}\int_Q\Exp\BNorm{\sum_{R\subseteq Q}
    \radem_R b_R(x)\ave{u}_R}{X}^p\ud x\Big)^{q/p}\\
  &\lesssim \sum_{k=0}^{\infty}\Big(\frac{1}{\abs{Q}}
    \int_Q (A2^k)^p\Exp\Babs{\sum_{R\in\mathscr{Q}_k}
      \radem_R b_R(x)}^p\ud x\Big)^{q/p}.
\end{split}\end{equation*}
For $k=0$, it is clear that
\begin{equation*}
  \frac{1}{\abs{Q}}\int_Q
  \Exp\Babs{\sum_{R\in\mathscr{Q}_0}
      \radem_R b_R(x)}^p\ud x
  \leq \Norm{b}{\Car^p(Q)}^p.
\end{equation*}
For $k\geq 1$, we have, using the definition and disjointness of
the cubes $S\in\mathscr{F}_{k-1}$,
\begin{equation*}\begin{split}
  \frac{1}{\abs{Q}}\int_Q
  \Exp\Babs{\sum_{R\in\mathscr{Q}_k}
      \radem_R b_R(x)}^p\ud x
  &\leq\sum_{S\in\mathscr{F}_{k-1}}\frac{1}{\abs{Q}}\int_S
  \Exp\Babs{\sum_{R\subseteq S}\radem_R b_R(x)}^p\ud x\\
&\leq   \frac{|\bigcup \limits _{S \in \mathscr{F}_{k-1}} S|}{\abs{Q}} 
    \Norm{b}{\Car^p(Q)}^p. 
\end{split}\end{equation*}

Since  \(\bigcup_{S \in \mathscr{F}_{k-1}} S \subseteq\{M_Ru>A\cdot 2^{k-1}\}\),  
it follows that
\begin{equation*}\begin{split}
  \Big(&\frac{1}{\abs{Q}}\int_Q\BNorm{\sum_{R\subseteq Q}
    \radem_R b_R(x)\ave{u}_R}{X}^p\Big)^{1/p}\\
  &\lesssim A\Norm{b}{\Car^p(Q)}\Big[1+\sum_{k=1}^{\infty}
    2^{kq}\Big(\frac{\abs{\{M_R u>A\cdot 2^{k-1}\}}}{\abs{Q}}
     \Big)^{q/p}\Big]^{1/q} \\
  &\lesssim A\Norm{b}{\Car^p(Q)}\Big[1
    +\int_0^{\infty} t^q\mu(\frac{M_R u}{A}>t)^{q/p}\frac{\ud
    t}{t}\Big]^{1/q},
\end{split}\end{equation*}
and the choice \(A=\Norm{M_R u}{L^{p,q}(\mu)}\) yields the asserted
bound  (using the fact that \(L^{p,p}(\mu) = L^{p}(\mu)\)) .
\end{proof}

\begin{theorem}\label{carlesonestimate}
Let $X$ be an RMF space, $1<p<\infty$, and $\epsilon>0$. Then
\begin{equation*}
  \Big(\int_{\R^n}\Exp\BNorm{\sum_{R\in\triangle}
    \radem_R b_R(x)\ave{u}_R}{X}^p\Big)^{1/p}
  \lesssim\Norm{b}{\Car^{p+\epsilon}(\R^n)}
  \Norm{u}{L^p(\R^n,X)},
\end{equation*}
for all \(u \in L^{p}(\R^{n};X)\).
We may take $\epsilon=0$ if $X$ has type $p$.
\end{theorem}

\begin{proof}
By standard considerations, it is easy to see that it suffices to
prove the estimate with a fixed dyadic cube $Q$ in place of $\R^n$
and $R\in\triangle$ replaced by $R\subseteq Q$. After dividing
this modified claim by $\abs{Q}^{1/p}$, the left-hand side becomes
identical with that in Lemma~\ref{lem:coreOfCarleson}, while the
right-hand side is
$\Norm{b}{\Car^{p+\epsilon}(Q)}\Norm{u}{L^p(\mu)}$.
 If \(X\) has type \(p\), the result with \(\epsilon=0\) thus follows from Lemma~\ref{lem:coreOfCarleson}. 
We now turn to the case where \(X\) has type \(t < p\). 

By the real method of interpolation, after
linearizing $M_R u$ in a standard manner, we have that
$\Norm{M_R u}{L^{p,q}(\mu)}\lesssim\Norm{u}{L^{p,q}(\mu,X)}$ for
the same $p$ and $1\leq q\leq\infty$. Thus
Lemma~\ref{lem:coreOfCarleson} shows that the bilinear map
\begin{equation}\label{eq:CarlesonMap}
  (b,u)\mapsto\sum_{R\subseteq Q}\radem_R b_R(\cdot)\ave{u}_R
\end{equation}
is bounded
\begin{equation}\label{eq:CarlesonBound}
  \Car^p(Q)\times L^{p,t}(\mu,X)\to L^p(\mu,\Rad(X))
\end{equation}
if $X$ has type $t\leq p$. 

If $X$ does not have type $p$, it nevertheless has type $1$.
For a small number $\epsilon>0$, we already know the following boundedness properties of the Carleson map \eqref{eq:CarlesonMap}:
\begin{equation}\label{eq:CarDeductions}\begin{split}
  &\Car^{p+\epsilon}(Q)\times L^{p+\epsilon,1}(\mu,X)\to
  L^{p+\epsilon}(\mu,\Rad(X)), \\
  &\Car^{p+\epsilon}(Q)\times L^{p-\epsilon,1}(\mu,X)\to
  L^{p-\epsilon}(\mu,\Rad(X)).
\end{split}\end{equation}
The second line  uses the embedding 
$\Car^{p+\epsilon}(Q)\subseteq\Car^{p-\epsilon}(Q)$. For a
fixed $b\in\Car^{p+\epsilon}(Q)$, the lines
\eqref{eq:CarDeductions} express the boundedness of the linear
operator $u\mapsto \sum_{R\subseteq Q}\radem_R
b_R(\cdot)\ave{u}_R$ between certain function spaces. Using the real interpolation results
\begin{equation*}\begin{split}
  &(L^{p+\epsilon,1}(\mu,X),L^{p-\epsilon,1}(\mu,X))_{\theta,p} = L^p(\mu,X) \\
  &(L^{p+\epsilon}(\mu,\Rad(X)),L^{p-\epsilon}(\mu,\Rad(X)))_{\theta,p} =L^p(\mu,\Rad(X))
\end{split}\end{equation*}
for appropriate \(\theta\in(0,1)\), we deduce the assertion.
\end{proof}


\section{Carleson measure estimate}\label{sec:conclusion}

In Section~\ref{sect:reduction}, we reduced the asserted inequality of Proposition~\ref{mainprop} to the estimation of the principal part~\eqref{eq:principal}. We have finally developed the required tools for dealing with this part in this final section.

Let us first see how to make use of the fact that we only need to consider $u\in\range(\Gamma)$. Since $\Gamma$ is a first-order constant-coefficient partial differential operator in \(L^p(\R^n;\C^N)\), it has the form \(\Gamma=\Gamma_0\nabla\), where \(\Gamma_0\in\bddlin(\C^{n};\C^N)\). Let us write \(W_{\Gamma}:=\range(\Gamma_0)\subseteq\C^N\), and let \(P_{\Gamma}\) be the orthogonal projection of \(\C^N\) onto this subspace. As before, we use the same symbol for its tensor extension to \(X^N\). Now, for \(u\in\range(\Gamma)\), we have
\begin{equation*}
  \gamma_{2^k}(x)A_{2^k} u(x)
  =\gamma_{2^k}(x)P_{\Gamma}A_{2^k} u(x)
  =\frac{\gamma_{2^k}(x)P_{\Gamma}}{\Norm{\gamma_{2^k}(x)P_{\Gamma}}{}}
    \Norm{\gamma_{2^k}(x)P_{\Gamma}}{}A_{2^k} u(x),
\end{equation*}
where we denote by \(\Norm{\gamma_{2^k}(x)P_{\Gamma}}{}\) the operator norm of \(\gamma_{2^k}(x)P_{\Gamma}\) in \(\bddlin(\C^N)\) (and let \(0/0:=0\)). Since the tensor extensions of the operators \(M\in\bddlin(\C^N)\) with \(\Norm{M}{}\leq 1\) are \(R\)-bounded on \(X^N\) (by writing out the matrix multiplications and using the contraction principle), it follows from Theorem \ref{carlesonestimate}
\begin{equation}\label{eq:ppToCarleson}\begin{split}
  \Exp &\BNorm{\sum_k\radem_k\gamma_{2^k}A_{2^k} u}{L^p(\R^n;X^N)}
  \lesssim\Exp\BNorm{\sum_k\radem_k\Norm{\gamma_{2^k}P_{\Gamma}}{}A_{2^k} u}{L^p(\R^n;X^N)} \\
  &=\Exp\BNorm{\sum_{Q\in\triangle}\radem_Q 1_Q\Norm{\gamma_{\ell(Q)}P_{\Gamma}}{}\ave{u}_Q}{L^p(\R^n;X^N)} \\
  &\lesssim\BNorm{\Big(1_Q\Norm{\gamma_{\ell(Q)}P_{\Gamma}}{}\Big)_{Q\in\triangle}}{\Car^{p+\epsilon}(\R^n)}
   \Norm{u}{L^p(\R^n;X^N)}.
\end{split}\end{equation}
Hence proving the asserted quadratic estimate in \(L^p(\R^n;X^N)\) is finally reduced to showing the finiteness of the \(\Car^{p+\epsilon}(\R^n)\)-norm above.

There are two peculiarities worth pointing out here. First, the space \(X\) has completely disappeared from this remaining estimate. Hence, the rest of the proof will be merely an \(L^p\) version, no longer Banach space valued, of the \(L^2\) estimates in \cite{AKM}.

Second, to get our desired \(L^p\) inequality, we are now required to prove an \(L^{p+\epsilon}\)-type estimate. This (and only this) is the reason why we formulated the main results --- Theorem~\ref{mainthm}, Corollary~\ref{maincor}, and Proposition~\ref{mainprop} --- for \(p\) in an open interval \((p_-,p_+)\), instead of just a single exponent \(p\). At this point it could seem that we only need openness at the upper end of the interval, but we also have to be able to repeat the reasoning in the dual case with the interval \((p_+',p_-')\).

The reader may also recall that the \(\epsilon\) could be avoided in \eqref{eq:ppToCarleson} if \(X\) has type \(p\). But to make the dual argument, we would also require that \(X^*\) has type \(p'\), and the only exponent for which this can be the case is \(p=2\). Moreover, if both \(X\) and \(X^*\) have type \(2\), then \(X\) is isomorphic to a Hilbert space, and so we are back to the classical situation. Thus we are able to recover the original \(L^2\) result in Hilbert spaces, but this is also the only situation, where we can work in a fixed \(L^p\) space. 

Now that we have assumed this extra \(\epsilon\), it is clear that completing the proof will only require the following. (Note also that \(R\)-bisectoriality of an operator \(T\otimes I_X\) in \(L^p(\R^n;X^N)\), where \(X\) is an arbitrary Banach space, implies \(R\)-bisectoriality of \(T\) in \(L^p(\R^n;\C^N)\) by restricting to a subspace.)

\begin{proposition}\label{carlesonmeasure}
Let \(1<p<\infty\), and let \(\Pi_B\) and \(\Pi_{B^*}\)  be perturbed Hodge--Dirac operators, which are \(R\)-bisectorial in \(L^p(\R^n;\C^N)\). Then
\begin{equation*}
  \Norm{\Big(1_Q \Norm{\gamma_{\ell(Q)}P_{\Gamma}
   }{\bddlin(\C^N)}
    \Big)_{Q\in\triangle}
     }{\Car^p(\R^n)}\lesssim 1.
\end{equation*}
\end{proposition}

The proof follows closely the Carleson measure estimate in Section~5 of \cite{AKM}, and hence we will skip some detail by simply asking the reader to repeat the relevant steps in~\cite{AKM}.

Denoting \(R_{Q} := (0,\ell(Q)]\times Q\), a reformulation of the claim is
\begin{equation*}\begin{split}
  &\Exp\BNorm{\sum_{k\in\Z}\radem_k 1_{R_Q}(2^k,\cdot)\gamma_{2^k}P_{\Gamma}}{L^p(\R^n;\bddlin(\C^N))} \\
  &\eqsim\BNorm{\Big(\sum_{k\in\Z}\Norm{1_{R_Q}(2^k,\cdot)\gamma_{2^k}P_{\Gamma}}{\bddlin(\C^N)}^2\Big)^{1/2}}{L^p(\R^n)}
  \lesssim \abs{Q}^{1/p}.
\end{split}\end{equation*}
The equivalence of the first and second form may be justified by Kahane's inequality and using the equivalent Hilbert--Schmidt norm on the finite-dimensional operator space \(\bddlin(\C^N)\).

Let us introduce the following subspace of \(\bddlin(\C^N)\), which contains our operators of interest \(\gamma_{2^k}(x)P_{\Gamma}\):
\begin{equation*}
  \mathscr{O}_{\Gamma}:=\{\nu\in\bddlin(\C^N):W_{\Gamma}^{\perp}\subseteq\kernel(\nu)\}
  =\{\nu\in\bddlin(\C^N):\nu=\nu P_{\Gamma}\}.
\end{equation*}
We set \(\sigma>0\) to be chosen later, and consider the cones
\begin{equation*}
  K_{\nu} = \Big\{\nu' \in \mathscr{O}_{\Gamma}\setminus\{0\} : \Norm{\frac{\nu'}{\Norm{\nu'}{}}-\nu}{}\leq \sigma\Big\},
\end{equation*}
where \(\nu\) belongs to a finite set \(\Lambda\) such that \(\bigcup_{\nu \in \Lambda}K_{\nu} = \mathscr{O}_{\Gamma}\setminus\{0\}\).
Writing
\begin{equation*}
  C_{\nu} := \{(t,x) \in (0,\infty) \times \R^{n} : \gamma_{t}(x)P_{\Gamma} \in K_{\nu}\},
\end{equation*}
we need to show that
\begin{equation*}
  \Exp\BNorm{\sum_{k \in \Z} \radem_{k} 1_{R_{Q} \cap C_{\nu}}(2^{k},.)\gamma_{2^{k}}P_{\Gamma}}{p} \lesssim |Q|^{1/p}
\end{equation*}
for each \(\nu\in\Lambda\). This in turns reduces to proving the following Proposition.

\begin{proposition}\label{lem:AKM59}
There exist \(\beta \in (0,1)\) and \(C>0\) which satisfy the following. For all \(Q \in \triangle\) and all \(\nu \in \bddlin(\C^n)\) with \(\Norm{\nu}{}=1\), there is a collection
\((Q_{j})_{j \in J}\) of disjoint dyadic subcubes of \(Q\) such that: denoting
\begin{equation}\label{eq:EQnu}
  E_{Q,\nu} := Q \setminus \bigcup_{j \in J}Q_{j},\qquad
  E^{*}_{Q,\nu} := R_{Q} \setminus \bigcup_{j \in J}R_{Q_{j}},
\end{equation}
there holds \(\abs{E_{Q,\nu}} > \beta \abs{Q}\) and
\begin{equation*}
  \Big(\Exp\BNorm{\sum_{k \in \Z} \radem_{k} 1_{E^{*}_{Q,\nu}\cap C_{\nu}}(2^{k},.)\gamma_{2^{k}}P_{\Gamma}}{p}^p
    \Big)^{1/p}\leq C\abs{Q}^{1/p}.
\end{equation*}
\end{proposition}

Indeed, assuming this is proven, we have for a fixed \(Q\in\triangle\)
\begin{equation*}
  \Exp \BNorm{\sum_{k \in \Z} \radem_{k}
     1_{R_{Q} \cap C_{\nu}}(2^{k},\cdot)\gamma_{2^{k}}P_{\Gamma}}{p}^p
  \leq C^p\abs{Q}+
  \sum_{j\in J}\Exp\BNorm{\sum_{k\in\Z}\radem_k
     1_{R_{Q_j}\cap C_{\nu}}(2^k,\cdot)\gamma_{2^k}P_{\Gamma}}{p}^p.
\end{equation*}
Now, applying Proposition \ref{lem:AKM59} for each of the \(Q_{j}\), and denoting by \((Q_{j,j'})_{j' \in J'}\) 
the corresponding sequence of subcubes of \(Q_{j}\), we have
\begin{equation*}\begin{split}
  \Exp &\BNorm{\sum_{k \in \Z} \radem_{k}
     1_{R_{Q} \cap C_{\nu}}(2^{k},\cdot)\gamma_{2^{k}}P_{\Gamma}}{p}^p \\
  &\leq C^p\abs{Q}+ C^{p}\sum_{j \in J}\abs{Q_{j}}+
 \sum_{j\in J}\sum_{j'\in J'}\Exp\BNorm{\sum_{k\in\Z}\radem_k
     1_{R_{Q_{j,j'}}\cap C_{\nu}}(2^k,\cdot)\gamma_{2^k}P_{\Gamma}}{p}^p.\\
  &    \leq C^p\abs{Q}(1+(1-\beta))+
 \sum_{j\in J}\sum_{j'\in J'}\Exp\BNorm{\sum_{k\in\Z}\radem_k
     1_{R_{Q_{j,j'}}\cap C_{\nu}}(2^k,\cdot)\gamma_{2^k}P_{\Gamma}}{p}^p.
\end{split}\end{equation*}

Reiterating this procedure leads to 
\begin{equation*}
 \Exp \BNorm{\sum_{k \in \Z} \radem_{k} 1_{R_{Q} \cap C_{\nu}}(2^{k},\cdot)\gamma_{2^{k}}P_{\Gamma}}{p}^p
  \leq C^{p}\abs{Q}\sum_{i=0}^{\infty} (1-\beta)^{i} = C^{p}\abs{Q}\beta^{-1}.
\end{equation*}

We now turn to the proof of Proposition~\ref{lem:AKM59}. Let us fix \(\nu\in\mathscr{O}_{\Gamma}\subseteq\bddlin(\C^N)\) of norm \(1\), and let \(w,\hat{w}\in\C^N\) also be of norm \(1\), and such that \(w=\nu^*(\hat{w})=P_{\Gamma}\nu^*(\hat{w})\). Hence \(w\in W_{\Gamma}\). We can now construct (as in \cite{AAH}, Lemma 4.10) the following kind of auxiliary functions for each \(Q\in\triangle\):
\begin{equation*}
  w_{Q} \in \range(\Gamma),\quad
  \supp w_Q\subseteq 3Q,\quad
  w_Q(x)\equiv w\ \forall x\in 2Q,\quad
  \Norm{w_Q}{\infty}\lesssim 1.
\end{equation*}
To do so, we take an affine function \(u_Q\) such that \(\Gamma u_Q\equiv w\) and \(\|1_{Q}u_Q\|_{\infty} \lesssim \ell(Q)\), and a smooth cutoff \(\eta_{Q}\) supported in \(3Q\) and equal to 1 on \(2Q\), with \(\|\nabla \eta_Q\|_{\infty} \lesssim \ell(Q)^{-1}\). Then we define 
\(w_{Q} = \Gamma(\eta_{Q}u_Q)\).

We now set \(f_{Q}^{w} := P_{\varepsilon\ell(Q)} ^{B}w_{Q}\). This satisfies
\begin{equation}\label{test1}
  \Norm{f_{Q}^{w}}{p} \lesssim \Norm{w_{Q}}{p} \lesssim \abs{Q}^{1/p},
\end{equation}
and, using the identity \(Q_s P_t=s/t\cdot Q_t P_s\), also
\begin{equation}\label{test2}\begin{split}
  \Exp\BNorm{\sum_{k \in \Z}\radem_k 1_{R_Q}(2^{k},.)Q_{2^{k}}^{B}f_{Q}^{w}}{p}
  \leq \sum_{k:2^{k} \leq \ell(Q)} \frac{2^{k}}{\varepsilon\ell(Q)}
   &\Norm{Q_{\varepsilon\ell(Q)}^{B} P_{2^{k}}^{B} w_{Q}}{p} \\
  &\lesssim \frac{\abs{Q}^{1/p}}{\varepsilon}.
\end{split}\end{equation}

Estimates \eqref{test1} and \eqref{test2} are our \(L^p\) versions of the first two assertions of \cite{AKM}, Lemma~5.10, and the remaining part of that Lemma is dealt with as follows. Note that we write simply \(\abs{\cdot}\) for the norm in \(\C^N\).

\begin{lemma}\label{lem:closeToAv}
For some $c$ depending only on $p$ as well as $P_t^B$, $Q_t^B$, and $\Gamma$, there holds
\begin{equation*}
  \Babs{\intav_Q f^w_Q\ud x-w}\leq c\varepsilon^{1/p'}.
\end{equation*}
\end{lemma}

\begin{proof}
Writing out the definitions,
\begin{equation}\label{eq:expressDiff}\begin{split}
  \intav_Q f^w_Q\ud x-w
  &=\intav_Q (P^B_{\varepsilon\ell(Q)}-I)w_Q\ud x \\
  &=\intav_Q -\varepsilon^2\ell(Q)^2\Gamma\Pi_B P^B_{\varepsilon\ell(Q)}w_Q\ud x,
\end{split}\end{equation}
where the last equality used the facts that $w_Q\in\range(\Gamma)$ and $\Pi_B^2=\Gamma\Pi_B$ on $\range(\Gamma)$. We next make use of the following estimate, which depends on the fact that $\Gamma$ is a first-order 
 differential operator with constant coefficients:
\begin{equation}\label{eq:Lp56}
  \Babs{\intav_Q\Gamma u\ud x}^p
  \lesssim\ell(Q)^{1-p}\Big(\intav_Q\abs{u}^p\ud x\Big)^{1/p'}
    \Big(\intav_Q\abs{\Gamma u}^p\ud x\Big)^{1/p}.
\end{equation}
This is the $L^p$ version of Lemma~5.6 in \cite{AKM}, and is proved by a simple modification of the $p=2$ case given there.

Using \eqref{eq:Lp56} in \eqref{eq:expressDiff}, we obtain
\begin{equation*}\begin{split}
  &\Babs{\intav_Q f^w_Q\ud x-w}^p\\
  &\lesssim \ell(Q)^{1-p}\Big(\intav\abs{\varepsilon\ell(Q)Q^B_{\varepsilon\ell(Q)}w_Q}^p\ud x\Big)^{1/p'}
    \Big(\intav\abs{(P^B_{\varepsilon\ell(Q)}-I)w_Q}^p\ud x\Big)^{1/p}\\
  &\lesssim \ell(Q)^{1-p}(\varepsilon\ell(Q))^{p/p'}\Big(\abs{Q}^{-1}\int\abs{w_Q}^p\ud x\Big)^{1/p'+1/p}
  \lesssim \varepsilon^{p-1}
\end{split}\end{equation*}
by the uniform $L^p$-boundedness of $P^B_t$ and $Q^B_t$, together with \eqref{test1}, and this completes the proof.
\end{proof}

\begin{lemma}
\label{lem:stoptime}
With \(\varepsilon=(2c)^{-p'}\), where \(c\) is as in Lemma~\ref{lem:closeToAv}, there exist \(\beta,c_1,c_2>0\) and for each \(Q\in\triangle\) a collection \((Q_j)_{j\in J}\) of disjoint dyadic subcubes such that, with the definitions \eqref{eq:EQnu}, there holds \(\abs{E_{Q,\nu}}>\beta\abs{Q}\) and
\begin{equation*}
  \Re(w,A_{2^k}f^w_Q(x))\geq c_1,\quad
  A_{2^k}\abs{f^w_Q}(x)\leq c_2,\qquad\text{if}\quad
  (2^k,x)\in E^*_{Q,\nu}.
\end{equation*}
\end{lemma}

\begin{proof}
With the given choice of \(\varepsilon\),  Lemma~\ref{lem:closeToAv} implies that
\begin{equation*}
  \Re\Big(w,\intav_{Q} f_{Q}^{w}\Big) \geq \frac{1}{2}.
\end{equation*}
The assertion follows from this together with \eqref{test1}, by a stopping time argument exactly as the corresponding result, Lemma~5.11, in \cite{AKM}.
\end{proof}

\begin{lemma}
\label{lem:test}
With \(\sigma:=\frac{c_1}{2c_2}\), there holds
\begin{equation*}
  \abs{\gamma_{2^k}(x)\big(A_t f^w_Q(x)\big)}\geq\frac{c_1}{2}\Norm{\gamma_{2^k}(x)P_{\Gamma}}{},\qquad
  (2^k,x)\in E^*_{Q,\nu}\cap C_{\nu}.
\end{equation*}
\end{lemma}

\begin{proof}
This is almost like \cite{AKM}, Lemma~5.12. By Lemma \ref{lem:stoptime},
\begin{equation*}
  \abs{\nu\big(A_{2^k}f^w_Q(x)\big)}
  \geq\Re\big(\hat{w},\nu\big(A_{2^k}f^w_Q(x)\big)\big)
  =\Re\big(w,A_{2^k}f^w_Q(x)\big)\geq c_1,
\end{equation*}
and then
\begin{equation*}\begin{split}
  &\abs{\frac{\gamma_{2^k}(x)P_{\Gamma}}{\Norm{\gamma_{2^k}(x)P_{\Gamma}}{}}
   \big(A_{2^k}f^w_Q(x)\big)} \\
  &\geq\abs{\nu\big(A_{2^k}f^w_Q(x)\big)}
  -\Norm{\frac{\gamma_{2^k}(x)P_{\Gamma}}{\Norm{\gamma_{2^k}(x)P_{\Gamma}}{}}-\nu}{}
   \abs{A_{2^k}f^w_Q(x)} \\
  &\geq c_1-\sigma c_2=c_1/2.
\end{split}\end{equation*}
Finally, recall that \(P_{\Gamma}\big(A_{2^k}f^w_Q(x)\big)=A_{2^k}f^w_Q(x)\), since \(f^w_Q\in\range(\Gamma)\), to complete the proof.
\end{proof}

\begin{proof}[of Proposition~\ref{lem:AKM59} and Proposition~\ref{carlesonmeasure}]
We make use of the Khintchine--Kahane inequalities (Proposition \ref{kk-inequalities}) and Lemma~\ref{lem:test} to the result:
\begin{equation*}\begin{split}
  \Big(\Exp &\BNorm{\sum_{k\in\Z}\radem_k 1_{R_Q\cap E^*_{Q,\nu}}(2^k,\cdot)\gamma_{2^k}P_{\Gamma}}{
    L^p(\R^n;\bddlin(\C^N))}^p\Big)^{1/p} \\
  &\eqsim\BNorm{\Big(\sum_{k\in\Z} 1_{R_Q\cap E^*_{Q,\nu}}(2^k,\cdot)\Norm{\gamma_{2^k}P_{\Gamma}}{}^2\Big)^{1/2}}{
    L^p(\R^n)} \\
  &\lesssim\Exp\BNorm{\sum_{k\in\Z}\radem_k 1_{R_Q}(2^k,\cdot)
    \gamma_{2^k}A_{2^k}f^w_Q}{L^p(\R^n;\C^N)} \\
  &\leq\Exp\BNorm{\sum_{k\in\Z}\radem_k \big(Q_{2^k}^B-\gamma_{2^k}A_{2^k}\big)f^w_Q}{L^p(\R^n;\C^N)}+\\
  &\qquad +\Exp\BNorm{\sum_{k\in\Z}\radem_k 1_{R_Q}(2^k,\cdot)
    Q_{2^k}^B f^w_Q}{L^p(\R^n;\C^N)}.
\end{split}\end{equation*}
Recalling again that \(f^w_Q\in\range(\Gamma)\), we may apply the reduction-to-principal part Theorem~\ref{reduction}, which shows that the first term on the right is dominated by \(\Norm{f^w_Q}{p}\lesssim\abs{Q}^{1/p}\). The second term is almost like the quadratic norm in Proposition~\ref{mainprop} which we started from but with the arbitrary \(X^N\)-valued function \(u\in\range(\Gamma)\) replaced by the deliberately constructed \(\C^N\)-valued test function \(f^w_Q\). And indeed the estimate for this test function, which we recorded in \eqref{test2}, is precisely what we need to complete the proof.
\end{proof}

\begin{proof}[of Proposition~\ref{mainprop} and Theorem~\ref{mainthm}] 
By Proposition~\ref{carlesonmeasure} and our analogue of Carleson's inequality (Theorem~\ref{carlesonestimate}) we have:
\begin{equation*}
  \Exp\BNorm{\sum_{k\in\Z}\radem_k
    \gamma_{2^{k}}A_{2^{k}}u}{L^p(\R^n;X^N)} \lesssim \|u\|_{L^p(\R^n;X^N)},
  \qquad \forall u\in\range(\Gamma).
\end{equation*}
Together with our quadratic \(T(1)\) Theorem~\ref{reduction}, this completes the proof of Proposition~\ref{mainprop}, and, as pointed out in Section~\ref{result}, of Theorem~\ref{mainthm}.
\end{proof}

\begin{remark}
Looking back at the structure of the entire proof, it may be interesting to note the difference in the two applications of Theorem~\ref{reduction}. In Section~\ref{sect:reduction}, it was used to replace \(Q_{2^k}^B\) in the desired estimate by its principal part \(\gamma_{2^k}A_{2^k}\), whereas right above we performed the reverse action. But of course other reductions took place at the same time: the first replacement allowed the application of Carleson's inequality, which reduced the original \(X^N\)-valued estimate to an \(\bddlin(\C^N)\)-valued one, while the second replacement made the further reduction to a \(\C^N\)-valued inequality for a test function. This strategy was already used in the case when \(X=\C\) in \cite{AKM}; thus the key point was not the reduction of \(X^N\) to \(\C^N\), but the reduction of \(u\) to \(f^w_Q\).
\end{remark}

\begin{appendix}

\section{$R$-bisectoriality of uniformly elliptic operators}\label{app:IplusK}

In this section we explain how the $R$-bisectoriality conditions in Theorem~\ref{mainthm} can, in some cases, be checked by a simple perturbation argument. Consider the differential operator $L=-\diverg A\nabla$, where the $\bddlin(\C^n)$-valued function $A(x)$ satisfies the uniform ellipticity (or accretivity) condition
\begin{equation}\label{eq:elliptic}
  \lambda\abs{\xi}^2\leq\Re\pair{A(x)\xi}{\xi},\qquad
  \abs{\pair{A(x)\xi}{\eta}}\leq \Lambda\abs{\xi}\abs{\eta}
\end{equation}
for all $x\in\R^n$ and $\xi,\eta\in\C^n$. This implies in particular that $x\mapsto A(x)$ and $x\mapsto A(x)^{-1}$ are in $L^{\infty}(\R^n;\bddlin(\C^n))$ with norms at most $\Lambda$ and $\lambda^{-1}$, respectively, as required to apply Corollary~\ref{maincor}. But the ellipticity \eqref{eq:elliptic} says more: as shown in \cite{mc90}, there exist constants $M,\delta>0$, depending only on $\lambda$ and $\Lambda$, such that $\Norm{MI-A(x)}{}\leq M-\delta$ for all $x\in\R^n$. Then $A=M(I+M^{-1}[A-MI])=:M(I+K)$, where the norm of $K$ in $L^{\infty}(\R^n;\bddlin(\C^n))$
is strictly smaller than~$1$. This obviously implies the same norm bound in $\bddlin(L^p(\R^n;\C^n))$. To be able to make this conclusion even in $\bddlin(L^p(\R^n;X^n))$, we need to use a special norm in the product space $X^n$. This is given by
\begin{equation}\label{eq:EuclNorm}
  \Norm{(x_i)_{i=1}^n}{X^n}
  :=\Big(\Exp\bnorm{\sum_{i=1}^n\gamma_i x_i}{X}^2\Big)^{1/2},
\end{equation}
where the $\gamma_i$ are independent standard Gaussian random variables. This is, of course, equivalent to any of the usual norms that one would use on $X^n$, and the equivalence constants may be chocen to depend on $n$ only. The crucial property of this norm is the following:

\begin{lemma}
Let $T\in\bddlin(\C^n)$ induce an operator in $\bddlin(X^n)$ in the natural way. If $X^n$ is equipped with the norm~\eqref{eq:EuclNorm}, then
\begin{equation*}
  \Norm{T}{\bddlin(X^n)}=\Norm{T}{\bddlin(\C^n)}.
\end{equation*}
\end{lemma}

\begin{proof}
The inequality $\geq$ is clear. The estimate $\leq$ follows from \cite{pisier:facto}, Proposition~3.7, once we observe that
\begin{equation*}\begin{split}
  \sum_{i=1}^n &\Babs{\pair{\sum_{j=1}^n t_{ij}x_j}{x^*}}^2
  =\norm{T(\pair{x_j}{x^*})_{j=1}^n}{\C^n}^2 \\
  &\leq\Norm{T}{\bddlin(\C^n)}^2\norm{(\pair{x_j}{x^*})_{j=1}^n}{\C^n}^2
  =\Norm{T}{\bddlin(\C^n)}^2\sum_{j=1}^n \abs{\pair{x_j}{x^*}}^2
\end{split}\end{equation*}
for all $x^*\in X^*$.
\end{proof}

We will now make use of the above observations but applied to $A^{-1}$ in place of $A$. Note that $A^{-1}$ also satifies the ellipticity condition \eqref{eq:elliptic}, possibly with different constants, as soon as $A$ does. Since the differential operators $L$ and $ML$ have the same mapping properties, we may assume without loss of generality that $M=1$. Thus the matrix-multiplication operator $A$ as in \eqref{eq:elliptic} may be assumed to have an inverse, which is a perturbation of the identity:
\begin{equation}\label{eq:IplusK}
  A^{-1}=I+K,\qquad\Norm{K}{\bddlin(L^p(\R^n;X^n))}\leq \Norm{K}{L^{\infty}(\R^n;\bddlin(\C^n))}<1.
\end{equation}

Hence, keeping the notation of Theorem~\ref{mainthm} and Corollary~\ref{maincor}, with \(A_{1}=I\) and \(A_{2}=A\),
\begin{equation}\label{eq:ellPiB}
  \Pi_B=\begin{pmatrix} 0 & -\diverg A \\ \nabla & 0 \end{pmatrix},\qquad
  \Pi_{B^*}=\begin{pmatrix} 0 & -\diverg \\ A\nabla & 0 \end{pmatrix}.
\end{equation}
and then
\begin{equation*}\begin{split}
  (I+&it\Pi_B)\begin{pmatrix} I & 0 \\ 0 & A^{-1} \end{pmatrix}
   =\begin{pmatrix} I & 0 \\ 0 & A^{-1} \end{pmatrix}(I+it\Pi_{B^*}) \\
  &=\begin{pmatrix} I & -it\diverg \\ it\nabla & A^{-1} \end{pmatrix}
   =\begin{pmatrix} I & -it\diverg \\ it\nabla & I \end{pmatrix}\left[
    I+\begin{pmatrix} I & -it\diverg \\ it\nabla & I \end{pmatrix}^{-1}
    \begin{pmatrix} 0 & 0 \\ 0 & K \end{pmatrix}\right] \\
  &=(I+it\Pi)\begin{pmatrix} I & it\diverg(I-t^2\nabla\diverg)^{-1}K \\ 
                             0 & I+(I-t^2\nabla\diverg)^{-1}K \end{pmatrix}.
\end{split}\end{equation*}

It follows that
\begin{equation}\label{eq:invertible}\begin{split}
  &(I+it\Pi_B)\text{ is invertible}\quad\Leftrightarrow\quad
  (I+it\Pi_{B^*})\text{ is invertible}\\
  &\qquad\Leftrightarrow\quad\left(I+(I-t^2\nabla\diverg)^{-1}K\right)\text{ is invertible},
\end{split}\end{equation}
and if this is the case, then
\begin{equation}\label{eq:ellResolvent}\begin{split}
  &\begin{pmatrix} I & 0 \\ 0 & A \end{pmatrix}R_t^B
   =R_t^{B^*}\begin{pmatrix}I & 0 \\ 0 & A \end{pmatrix}\\
  &=\begin{pmatrix} I & it\diverg(I-t^2\nabla\diverg)^{-1}K \\ 0 & I \end{pmatrix}
   \begin{pmatrix} I & 0 \\ 0 & [I+(I-t^2\nabla\diverg)^{-1}K]^{-1} \end{pmatrix}
   R_t
\end{split}\end{equation}
where, we recall, \(R_t^B=(I+it\Pi_B)^{-1}\), \(R_t=(I+it\Pi)^{-1}\).

We can now conclude the following:

\begin{proposition}
Let $X$ be a UMD space, $1<p<\infty$, and $A\in L^{\infty}(\R^n;\bddlin(\C^n))$ satisfy~\eqref{eq:IplusK}. Then the operators $\Pi_B$ and $\Pi_{B^*}$ in \eqref{eq:ellPiB} are $R$-bisectorial in the space $L^p(\R^n;X^{n+1})$ provided that $I+(I-t^2\nabla\diverg)^{-1}K$ is invertible in $L^p(\R^n;X^n)$ for all $t>0$, and
\begin{equation*}
  \{[I+(I-t^2\nabla\diverg)^{-1}K]^{-1}\}_{t>0}\quad\text{is $R$-bounded in}\quad
  L^p(\R^n;X^n).
\end{equation*}
Hence, if the above condition is valid in an interval \((p-\varepsilon,p+\varepsilon)\), then $\Pi_B$ and $\Pi_{B^*}$ have an \(H^{\infty}\) functional calculus
 in $L^p(\R^n;X^{n+1})$, \(L\) has an \(H^{\infty}\) calculus  in $L^p(\R^n;X)$, and \(L\) satisfies Kato's square root estimates \(\|\sqrt{L}u\|_{p} \eqsim \|\nabla u\|_{p}\) for all \(u \in L^{p}(\R^{n};X)\).
\end{proposition}

\begin{proof}
We have already seen that the invertibility condition is both necessary and sufficient for the existence of the resolvents appearing in the definition of bisectoriality. If $X$ is a UMD space, then the unperturbed operator $\Pi$ is $R$-bisectorial, and  moreover the family of operators
\begin{equation*}
  \{it\diverg(I-t^2\nabla\diverg)^{-1}\}_{t>0}
  =\{it(I-t^2\Delta)^{-1}\diverg\}_{t>0}
\end{equation*}
is $R$-bounded from $L^p(\R^n;X)$ to $L^p(\R^n;X^n)$ (by Proposition~\ref{RbddMult}, since these are Fourier multiplier operators whose symbols have uniformly bounded variation). From \eqref{eq:ellResolvent}, and the fact that products of $R$-bounded sets remain $R$-bounded, we conclude the first assertion. The second is a consequence of Theorem \ref{mainthm} and Corollary \ref{maincor}.
\end{proof}

\begin{remark}\label{rem:Neumann}
If $n=1$, then the equivalent invertibility conditions in \eqref{eq:invertible} are always satisfied in $L^p(\R;X^2)$ resp. $L^p(\R;X)$, for all Banach spaces $X$ and all $p\in\left[1,\infty\right]$. In fact, in this case $(I-t^2\nabla\diverg)^{-1}=(I-t^2\Delta)^{-1}=\mathscr{P}_t$ is the convolution operator with kernel $(2t)^{-1}e^{-\abs{x}/t}$. This operator contracts all $L^p$ spaces, and hence $I+\mathscr{P}_tK$ has a bounded inverse represented by the convergent Neumann series
\begin{equation}\label{eq:Neumann}
  (I+\mathscr{P}_tK)^{-1}
  =\sum_{k=0}^{\infty}(-\mathscr{P}_tK)^k,
\end{equation}
since the operator norm of $K$ satisfies $\Norm{K}{}<1$.
\end{remark}

\begin{corollary}
Let $X$ be a UMD function lattice. Let $A\in L^{\infty}(\R^n;\C)$ satisfy~\eqref{eq:IplusK}. Then the operators $\Pi_B$ and $\Pi_{B^*}$ in \eqref{eq:ellPiB} are $R$-bisectorial in $L^p(\R^n;X^2)$ for all $p\in\left]1,\infty\right[$, and hence $L=-d/dx\,A(x)\,d/dx$  has an \(H^{\infty}\) calculus and satisfies the Kato's square root estimates in \(L^{p}(\R^{n};X)\), for all $p\in\left]1,\infty\right[$.
\end{corollary}

\begin{proof}
By Remark~\ref{rem:Neumann} and \eqref{eq:invertible}, we already know that the required resolvents exist. To prove the $R$-boundedness of $(I+\mathscr{P}_tK)^{-1}$, it suffices to show that the $R$-bounds of the terms in the Neumann series \eqref{eq:Neumann} converge. Let us investigate the $k$th term. Our aim is to show that
\begin{equation}\label{eq:RboundToShow}
  \Exp\BNorm{\sum_j\radem_j(\mathscr{P}_{t_j}K)^k u_j}{L^p(\R;X)}
  \lesssim\Norm{K}{\infty}^k\Exp\BNorm{\sum_j\radem_j u_j}{L^p(\R;X)},
\end{equation}
since this would allow us to sum up the series in $k$. Since $X$ is a function lattice with finite cotype, \eqref{eq:RboundToShow} is equivalent to the quadratic estimate
\begin{equation}\label{eq:quadraticToShow}
  \BNorm{\Big(\sum_j\abs{(\mathscr{P}_{t_j}K)^k u_j}^2\Big)^{1/2}}{p}
  \lesssim\Norm{K}{\infty}^k\BNorm{\Big(\sum_j\abs{u_j}^2\Big)^{1/2}}{p}.
\end{equation}

Let us denote the convolution kernel of $\mathscr{P}_t$ by $p_t(x):=(2t)^{-1}e^{-\abs{x}/t}$. The positivity of this function is of essential importance in what follows. Now
\begin{equation*}\begin{split}
  &\abs{(\mathscr{P}_t K)^k u(x)}\\
  &=\abs{\int\cdots\int p_t(x-y_1)K(y_1)\cdots
     p_t(y_{k-1}-y_k)K(y_k)u(y_k)\ud y_1\cdots\ud y_k} \\
  &\leq\int\cdots\int p_t(x-y_1)\abs{K(y_1)}\cdots
     p_t(y_{k-1}-y_k)\abs{K(y_k)}\abs{u(y_k)}\ud y_1\cdots\ud y_k \\
  &\leq\Norm{K}{\infty}^k\int\cdots\int p_t(x-y_1)\cdots
     p_t(y_{k-1}-y_k)\abs{u(y_k)}\ud y_1\cdots\ud y_k \\
  &=\Norm{K}{\infty}^k\mathscr{P}_t^k\abs{u}(x).
\end{split}\end{equation*}
Hence we have
\begin{equation*}
  \BNorm{\Big(\sum_j\abs{(\mathscr{P}_{t_j}K)^k u_j}^2\Big)^{1/2}}{p}
  \leq\Norm{K}{\infty}^k\BNorm{\Big(\sum_j(\mathscr{P}_{t_j})^k \abs{u_j})^2\Big)^{1/2}}{p}.
\end{equation*}
The right-hand side above is dominated by the right-hand side of \eqref{eq:quadraticToShow}, with the implied constant independent of $k$, since the two-parameter family of operators $\{\mathscr{P}_t^k:t>0,k\in\Z_+\}$ is $R$-bounded in $L^p(\R;X)$. In fact, these are Fourier multiplier operators with symbols $(1+t^2\abs{\xi}^2)^{-k}$, and one readily checks that they all have uniformly bounded variation, so that we may apply Proposition~\ref{RbddMult}.

This completes the proof of the $R$-bisectoriality. The final claim concerning the functional calculus and the Kato estimates is just an application of Theorem~\ref{mainthm} and Corollary~\ref{maincor}.
\end{proof}

Note that the $R$-boundedness of $\{\mathscr{P}_t^k:t>0,k\in\Z_+\}$, which played a r\^ole above, is still true in arbitrary UMD spaces; however, without the possibility of replacing the randomized norms by quadratic ones, there does not seem to be a way of extracting the $K$'s out of the operator product $(\mathscr{P}_tK)^k$. In the noncommutative \(L^p\) spaces, there are also versions of square functions available, but the proof above does not apply, since the modulus \(\abs{\cdot}\) does not satisfy the triangle inequality.

In general, the Neumann series argument shows that \(\Pi_{B^{*}}\) and \(\Pi_{B}\) are bisectorial provided the set
\[\{(I-t^{2}\nabla\diverg)^{-1}K; t\in \R\}\]
is R-bounded with constant \(c<1\).
If \(X\) is a Hilbert space, and \(p=2\), the R-bounds are just uniform bounds and thus \(c \leq \|K\|_{\mathcal{L}(L^{p}(\R^{n};X))} < 1\).
This gives back the solution of the Kato problem from \cite{AHLMT}.
Still in the Hilbertian situation, this also implies that, given a perturbation, there exists an open interval \((p^{A}_{-},p^{A}_{+}) \subset (1,\infty)\) containing \(2\) such that \eqref{eq:PiBRsect} holds.
This coincides with results from \cite{auscher}.
Computing the precise values of \(p^{A}_{-}\) and \(p^{A}_{+}\) seems, unfortunately, to be difficult.


\section{Carleson's inequality and paraproducts}\label{app:paraproducts}

Let us point out some consequences of Theorem~\ref{carlesonestimate} concerning vector-valued paraproducts
\begin{equation*}
  P(f,u):=\sum_{Q\in\triangle}
  \sum_{\eta}
  \frac{\pair{f}{h_Q^{\eta}}\ave{u}_Q}{\abs{Q}}h_Q^{\eta}.
\end{equation*}
These operators play the important r\^ole of principal parts of Calder\'on--Zygmund operators in the \(T(1)\) and \(T(b)\) theorems. Versions of these theorems in UMD spaces have been proved in \cite{Figiel:90,Hytonen:Tb,HytWeis}

The basic mapping property in the scalar case \(X=\C\) is
\begin{equation}\label{eq:paraBasic}
  \Norm{P(f,u)}{L^p(\R^n)}\lesssim\Norm{f}{BMO(\R^n)}\Norm{u}{L^p(\R^n)},
  \qquad 1<p<\infty.
\end{equation}
This reduces to the classical Carleson inequality for $p=2$, and may be extrapolated to the whole range $1<p<\infty$ by standard Calder\'on--Zygmund techniques. Alternatively, one may establish the $L^2$ estimate in all weighted spaces $L^2(\R^n,w(x)\ud x)$ for $w$ in the Muckenhoupt $A_2$-class, with uniform dependence on the $A_2$-constant, and invoke the weighted extrapolation theorem of Rubio de Francia to deduce the corresponding $L^p$-estimates (cf.~\cite{KP} for this approach). Figiel \cite{Figiel:90} has shown (based on an intermediate estimate \cite{FigWoj}, which he attributes to Bourgain) that one may replace $L^p(\R^n)$ by $L^p(\R^n;X)$ in \eqref{eq:paraBasic} provided that $X$ is a UMD space. His proof employs interpolation between $(H^1,L^1)$ and $(L^{\infty},BMO)$ type estimates. Thus in all these arguments, the $L^p$-inequalities in \eqref{eq:paraBasic} when $p\neq 2$ are reached somewhat indirectly.

We next provide an alternative approach to the Bourgain--Figiel result based on Theorem~\ref{carlesonestimate} (and hence under the additional assumption of the RMF property). This also gives an apparently new ``\(L^p\) proof'' of the classical estimate \eqref{eq:paraBasic}. While the proof of Theorem~\ref{carlesonestimate} was not completely interpolation-free, either, one should note that getting the $L^p$ estimate for a given $p$ only involved interpolation between spaces ``in the proximity'' of $L^p$, in contrast to the ``far away'' end-point spaces in the classical arguments. The proof below will show that the problem of the extra $\epsilon$ disappears in this specific situation, thanks to the John--Nirenberg inequality. 

\begin{corollary}\label{cor:paraproducts}
Let $X$ be a UMD space with RMF, and $1<p<\infty$. Then
\begin{equation*}
  \Norm{P(f,u)}{L^p(\R^n;X)}
  \lesssim\Norm{f}{BMO(\R^n)}\Norm{u}{L^p(\R^n;X)}.
\end{equation*}
\end{corollary}

\begin{proof}
We have the following chain of estimates, where we write simply \(\Norm{\cdot}{p}\) for the norm of \(L^p(\R^n;X)\):
\begin{equation*}\begin{split}
  &\Norm{P(f,u)}{p} \\
  &\lesssim\Big(\int_{\R^n}\Exp\BNorm{
    \sum_{Q,\eta}\radem_Q^{\eta} 
    \frac{\pair{f}{h_Q^{\eta}}h_Q^{\eta}(x)
   }{\abs{Q}}\ave{u}_Q}{X}^p\ud x\Big)^{1/p}\\
  &\lesssim\sum_{\eta}\sup_{S\in\triangle}\Big(\frac{1}{\abs{S}}
   \int_S\Exp\Babs{\sum_{Q\subseteq S}\radem_Q
    \frac{\pair{f}{h_Q^{\eta}}h_Q^{\eta}(x)}{\abs{Q}}}^{p+\epsilon}
     \ud x\Big)^{1/(p+\epsilon)}
  \Norm{u}{p} \\
  &\lesssim\sup_{S\in\triangle}\Big(\frac{1}{\abs{S}}
   \int_S\Babs{\sum_{Q\subseteq S}\sum_{\eta}
    \frac{\pair{f}{h_Q^{\eta}}h_Q^{\eta}(x)}{\abs{Q}}}^{p+\epsilon}
    \ud x\Big)^{1/(p+\epsilon)}
  \Norm{u}{p} \\
  &=\sup_{S\in\triangle}\Big(\frac{1}{\abs{S}}
   \int_S\abs{f(x)-\ave{f}_S}^{p+\epsilon}
    \ud x\Big)^{1/(p+\epsilon)}
  \Norm{u}{p} \\
  &\lesssim\Norm{f}{BMO}\Norm{u}{p}.
\end{split}\end{equation*}
The first estimate employed the UMD property of $X$, the second used Theorem~\ref{carlesonestimate}, the third the UMD property of $\C$, and the final one the John--Nirenberg inequality.
\end{proof}

It is also possible to reverse the r\^oles of scalar and vector-valued functions in Theorem~\ref{carlesonestimate} and then in Corollary~\ref{cor:paraproducts}. We leave the straightforward verification of the details to the reader, and only record the result. The RMF property does not enter this time, because the maximal function estimate is now required for a scalar-valued function.

\begin{corollary}
Let $X$ be a UMD space, and $1<p<\infty$. Then
\begin{equation*}
  \Norm{P(f,u)}{L^p(\R^n;X)}
  \lesssim\Norm{f}{BMO(\R^n;X)}\Norm{u}{L^p(\R^n)}.
\end{equation*}
\end{corollary}

\section{The space $\ell^1$ does not have RMF}\label{app:l1}
As mentioned in Section \ref{maximal}, we do not yet understand how the RMF property relates to other properties of Banach spaces, and in particular to the UMD property. In this Appendix we show that it is, however, a nontrivial property by proving that \(\ell_{1}\) does not enjoy RMF.

Let \(n \in \N\), and $u(x)=e_k$ for
$x\in\left[(k-1)2^{-n},k2^{-n}\right)$ for $k=1,2,\ldots,2^n$. Then $\Norm{u}{L^p(\R^1,\ell^1)}
=1$ for all $p\in[1,\infty]$. For
\(x\in\left[0,2^{-n}\right)\), we have
\begin{equation*}
  A_{2^{-n+j}}u(x)=\frac{1}{2^j}\sum_{k=1}^{2^j}e_k,\qquad
  j=0,1,\ldots,n.
\end{equation*}
For other $x\in[0,1)$, we have similar results with a permuted basis $e_{\pi(k)}$ in place of $e_k$.

Let \(n=2^m\), and consider, given a sequence \(\alpha = (\alpha_{i})_{\i \in \N} \subset \R\) to be chosen later, the sequence \(\lambda\) given by \(\lambda_{2^i}=\alpha_i\), \(i=1,\ldots,m\), and \(\lambda_j=0\) otherwise. Then for \(0<x<2^{-n}\),
\begin{equation*}\begin{split}
  \Exp\BNorm{\sum_{j=0}^n\radem_j A_{2^{-n+j}}u(x)\lambda_j}{\ell^1}
  &=\Exp\BNorm{\sum_{i=1}^m\radem_i\frac{1}{2^{2^i}}\sum_{k=1}^{2^{2^i}}e_k\alpha_i}{\ell^1} \\
  &\geq\Exp\BNorm{\sum_{i=1}^m\radem_i\frac{1}{2^{2^i}}\sum_{k=2^{2^{i-1}}+1}^{2^{2^i}}e_k\alpha_i}{\ell^1}
   -\sum_{i=1}^m\frac{1}{2^{2^i}}2^{2^{i-1}}\abs{\alpha_i}\\
  &=\sum_{i=1}^m\frac{2^{2^i}-2^{2^{i-1}}}{2^{2^i}}\abs{\alpha_i}-\sum_{i=1}^m 2^{-2^{i-1}}\abs{\alpha_i}\\
  &\gtrsim\Norm{\alpha}{\ell^1}-\Norm{\alpha}{\ell^{\infty}}.
\end{split}\end{equation*}
Choosing, say, \(\alpha_i=(i+1)^{-1}\), we find that
\begin{equation*}
  M_Ru(x)\gtrsim\log m\gtrsim\log\log n
\end{equation*}
for all \(x\in[0,2^{-n})\), and by the permutation symmetry of the standard basis, for all \(x\in[0,1)\). This shows that \(\Norm{M_Ru}{L^p(\R^1)}\gtrsim\log\log n\). Since the same construction can be repeated with arbitrarily large \(n\), we see that no \(L^p\) bound can hold for \(M_R\) in \(\ell^1\).\\

\end{appendix}

\subsubsection*{Acknowledgments.}
Alan McIntosh and Pierre Portal would like to thank the Centre for Mathematics and its Applications at the Australian National University, and the Australian Research Council for their support. Tuomas Hyt\"onen gratefully acknowledges the support of the Finnish Academy of Science and Letters (Vilho, Yrj\"o and Kalle V\"ais\"al\"a Foundation), and the Academy of Finland (project 114374 ``Vector-valued singular integrals'').

\end{document}